\documentclass[12pt,reqno,a4paper]{amsart}

\usepackage{color}

\usepackage[margin=2cm]{geometry}

\newtheorem{theorem}{Theorem}[section]

\newtheorem{corollary}[theorem]{Corollary}
\newtheorem{proposition}[theorem]{Proposition}

{\theoremstyle{definition}

}

\newcommand{\R}{\mathbb{R}}

\newcommand{\T}{\theta}

\newcommand{\CO}{\mathcal{O}}

\DeclareMathOperator{\Span}{Span}
\DeclareMathOperator{\sgn}{sgn}
\DeclareMathOperator{\ee}{e}

\begin{document}
\title[Local versus global]{Local and global analysis of the displacement map for some near integrable systems}

\author[F. Braun]{Francisco Braun}
\address{Departamento de Matem\'{a}tica, Universidade Federal de S\~ao Carlos, 13565--905 S\~ao Carlos, S\~ao Paulo, Brazil}
\email{franciscobraun@ufscar.br} 

\author[L.P.C. da Cruz]{Leonardo Pereira Costa da Cruz}
\address{Instituto de Ci\^encias Matem\'aticas e Computa{\c{c}}{\~a}o, Universidade de São Paulo,  13566--590 S\~ao Carlos, S\~ao Paulo, Brazil}
\email{leonardocruz@icmc.usp.br}

\author[J. Torregrosa]{Joan Torregrosa}
\address{Departament de Matem\`{a}tiques, Universitat Aut\`{o}noma de Barcelona and Centre de Recerca Matem\`atica, 08193 Be\-lla\-ter\-ra, Barcelona, Spain}
\email{joan.torregrosa@uab.cat}



\date{\today}

\begin{abstract}
In this paper, we introduce an alternative method for applying averaging theory of orders $1$ and $2$ in the plane. 
This is done by combining Taylor expansions of the displacement map with the integral form of the Poincaré--Poyntriagin--Melnikov function. 
It is known that, to obtain results of order $2$ with averaging theory, the first-order averaging function should be identically zero. 
However, when working with Taylor expansions of the $i$th-order averaging function, we usually cannot guarantee it is identically zero. 
We prove that the vanishing of certain coefficients of the Taylor series of the first-order averaging function ensures it is identically zero. 
We present our reasoning in several concrete examples: a quadratic Lotka--Volterra system, a quadratic Hamiltonian system, the entire family of quadratic isochronous differential systems, and a cubic system. 
For the latter, we also show that a previous analysis contained in the literature is not correct. 
In none of the examples is it necessary to precisely calculate the averaging functions. 
\end{abstract}

\maketitle

\section{Introduction} \label{se:1}
Our object of study is a planar polynomial system of degree $n$, that is, a system of ordinary differential equations given by  
\begin{equation}\label{eq:0}  
	\left(\begin{array}{c}  
		\dot{x} \\  
		\dot{y}  
	\end{array}\right) = Z(x,y)= \left(P(x,y), Q(x,y)\right)^t,  
\end{equation}  
where $P$ and $Q$ are polynomials of degree $n$. The maximum number of limit cycles such a system can have is called the Hilbert number, denoted by $H(n)$. Finding $H(n)$ is a classical problem, at least since Hilbert posed it as part of Problem $16$ of his famous list of open problems for the $20$th century (see more details in \cite{Ily2003}). Although it is well known that $H(1) = 0$, the value of $H(n)$ remains open for all $n \geq 2$. The finiteness of the number of limit cycles for each vector field was established by Écalle \cite{Eca1992} and Ilyashenko \cite{Ily1991}. 
However, very recently, Yeung \cite{Yeu2025} has shown that the main technique used in Ilyashenko’s proof contains a mistake. Over the past several years, lower bounds for $H(n)$ have been found. 
On the other hand, upper bounds for $H(n)$ are known only for families of polynomial systems of degree $n$. Hence, the main problem is to establish the existence of a uniform upper bound.

In the study of the search for possible values of $H(n)$, several different ideas and techniques have been developed over the years. One of them is to make polynomial perturbations of centers in order to break the continuum of periodic orbits while preserving some isolated ones. This is the idea we employ in this work. We assume that the origin of coordinates is a non-degenerate center of the polynomial vector field $Z$ of \eqref{eq:0}. 
A polynomial perturbation of degree $n$ of $Z$ up to order $\ell$ is the vector field  
\begin{equation}\label{eq:2}  
	Z_{\ell,\varepsilon} (x,y) = Z(x,y) + \sum_{i=1}^\ell \varepsilon^i Z_i(x,y),  
\end{equation}  
where $\varepsilon > 0$, and $n,\ell$ are positive integers, and $Z_i(x,y) = \left(P_i(x,y), Q_i(x,y)\right)^t$, $i = 1,2$, are polynomial vector fields of degree $n$ without constant terms. This is not restrictive because, as this equilibrium point remains after perturbation, translating if necessary, we can assume it is fixed at the origin.

We also assume $Z(x,y) = R(x,y)(-\partial H/\partial y, \partial H/\partial x)$,  
with $H:\R^2 \rightarrow \R$ having an isolated minimum at the origin, such that $H(0,0) = 0$, and $R(x,y)$ is a suitable function such that $Z(x,y)$ is polynomial. Then we parameterize the center by the energy levels $H(x,y) = h$, and the connected component corresponding to the closed orbit passing through the point $(r,0)$ is denoted by $\Gamma(r)$. Clearly, $H$ is a first integral of the unperturbed system.

As usual, we write the displacement map of system \eqref{eq:2} as  
\[
\mathcal{M}(r,\varepsilon) = \sum_{i = 1}^\ell \varepsilon^i \mathcal{M}_i(r) + \CO(\varepsilon^{\ell+1}).  
\]  
The expression of $\mathcal{M}_i(r)$ is a generalized Abelian integral, and it is known as the $i$th Melnikov function, or $i$th-order averaging function. 
In fact, the first-order averaging function  
\begin{equation}\label{1o}  
\mathcal{M}_1(r) = \int_{\Gamma(r)} (Q_1\,dx - P_1\,dy)/R  
\end{equation}  
is also known as the Poincaré--Poyntriagin--Melnikov function. Clearly, the simple zeros of $\mathcal{M}_1(r)$ will correspond, by the Implicit Function Theorem, to limit cycles of \eqref{eq:2}. When $\mathcal{M}_1(r) \equiv 0$, the simple zeros of $\mathcal{M}_2(r)$ will correspond to limit cycles, etc. For more details, we refer the reader to the classical references \cite{Per2001, Rou1998}.

This method started with the classical works of Lagrange and Laplace, who provided an intuitive justification of the mechanism. The first formalization of this procedure was given by Fatou \cite{FatouP1928}. This approach has been and continues to be used intensively when trying to find limit cycles by perturbing a periodic annulus.  
However, applying it is often very difficult, as it usually depends on integrals that have no explicit solution. The few cases in which one can find $\mathcal{M}_1$ explicitly are those in which it is possible to obtain the explicit parameterization of $\Gamma(r)$. For instance, the linear canonical center, which has a first integral of the form $H(x,y) = x^2 + y^2$, and a parameterization given by the classic polar change of coordinates, or the potential case with a first integral of the form $H(x,y) = y^2 + V(x)$, where $V(x)$ is a polynomial function. Further, in a few cases, it is possible to use Picard-Fuchs theory, see again \cite{Rou1998}. 
In short, the cases when we are able to obtain the Abelian integral or the first-order averaging function are rare.

Anyway, when one is interested only in the zeros of the averaging function close to $0$, as is the case when analyzing the birth of limit cycles after perturbing a center around an equilibrium point, the most common technique is to expand this function in a Taylor series around $0$. 
As we will see further, in order to simplify the technique, we compute the Taylor series of \eqref{1o} at the origin without knowing it explicitly. 
But it is clear we now face a new problem: this approach works smoothly for the first-order, but how can we now reach higher-order averaging functions and use them to get limit cycles? Because to use the second-order, we need the first-order averaging function to vanish identically. 
But we have only a $j$-jet of $\mathcal{M}_1$, and it is clear that in order to prove that it is identically zero, we need all the coefficients to vanish. 
That is, the vanishing of the first $j$ terms of this series does not imply the whole series is identically zero. 
Actually, this is theoretically true in our situation: since, as we will see later, the coefficients are algebraic expressions in terms of the perturbative coefficients, the Hilbert basis theorem guarantees that only a finite number of them being equal to zero implies that all of them are equal to zero. 
The question is: how many of them do we have to annihilate in each case? And how can we be sure that if $j$ of them vanish, all the others are zero as well?

In this paper we will be able to prove, for some quadratic families, that $\mathcal{M}_1$ is identically zero after annihilating the first coefficients of its Taylor expansion. 
Moreover, we will be able to predict the number of coefficients that need to be annulled in order to achieve this, without the need for explicit calculations of $\mathcal{M}_i$ expressions. 

The idea is to illustrate the power of our approach, to be explained in details in the next sections, perturbing some polynomial families of low degree:
\begin{equation}\label{eq:2_1} 
\begin{aligned} 
LV: & \left\{\begin{aligned}
\dot x & = -y (1+x),\\ 
\dot y & = x(1+y),
\end{aligned}\right. 
& \mathcal{H}: & \left\{\begin{aligned}
\dot x & =-y -\frac{\alpha}{2} x^2-\beta x y-\frac{3 \gamma}{2} y^2,\\
\dot y & =x +\frac{3\delta}{2} x^2 + \alpha x y+ \frac{\beta}{2} y^2, 
\end{aligned}\right. \\
CR_1: & \left\{\begin{aligned}
\dot x &= - y(1 - 2 \alpha x - 2 x^2),\\
\dot y &=x+ \alpha (y^2-x^2)+2 xy^2,
\end{aligned}\right. 
& & 
\end{aligned}
\end{equation}
where $\alpha, \beta, \gamma, \delta$ are real values, and 
\begin{equation}\label{eq:2_2}
\begin{aligned}
& S_1: \left\{\begin{aligned}
\dot x & = -y + x^2 - y^2, \\
\dot y & =x(1+2y), 
\end{aligned}\right. 
& \phantom{meio} 
& S_2: \left\{\begin{aligned}
\dot x & =-y+x^2, \\
\dot y & =x(1+y), 
\end{aligned}\right. 
\\
& S_3: \left\{\begin{aligned}
\dot{x} & = -y-\frac{4}{3} x^2,\\
\dot{y} & = x \left(1- \frac{16}{3} y\right),
\end{aligned}\right.
& \phantom{meio} 
&  S_4: \left\{\begin{aligned}
\dot x &=-y+\frac{16}{3} x^2-\frac{4}{3} y^2,\\
\dot y &=x\left(1 + \frac{8}{3} y\right).
\end{aligned}\right.
\end{aligned}
\end{equation} 
It might be clear that this approach can be further applied for other systems, smooth or piecewise situations, or higher dimensions. 
The three systems in \eqref{eq:2_1} are a Lotka-Volterra,  Hamiltonian, and a cubic reversible of the Chavarriga and Sabatini compilation \cite{ChaSab1999}, respectively.
The systems in \eqref{eq:2_2} are the isochronous quadratic systems identified by Loud \cite{Lou64}, see also \cite{MarRouTon1995}, written after suitable linear changes of variables. 

These examples were chosen because of the particularities of each of them. First, to apply the averaging method to the center $LV$, we need to deal with an integrating factor, which makes the calculation of the integrals more complicated. In the Hamiltonian $\mathcal{H}$ case, we do not have the explicit parameterization of the solution of the center, making the explicit calculations related to the averaging method very difficult. The reversible case $CR_1$ was studied for $\alpha = 0$ in \cite{LiZha2014} (the system in \cite{LiZha2014} is the system $CR_1$ after a linear change of variables). We studied it for any $\alpha$. Finally, perturbations of the centers $S_1$, $S_2$, and $S_3$ frequently appear in the literature, so we are able to compare our approach with them. For the system $S_4$, there are few studies—see one in \cite{ShiJiaXia2001}—mainly because it does not have a birational linearization.  
For this reason, we will focus on it.

Before presenting the main results of the paper, we will introduce some notation.  
We begin by writing the series of the \emph{$i$th-order averaging function}, $\mathcal{M}_i$ for $i=1,\ldots$, in the neighborhood of $0$ as  
\begin{equation}\label{Mism}
	\mathcal{M}_i(r) = M_{i}^{[j]}(r) + \CO(r^{j+1}), \ \ \text{where} \ \ 
	M_{i}^{[j]}(r) = \sum_{k=1}^j m_{i,k} r^k. 
\end{equation}
That is, $M_{i}^{[j]}$ is the $j$-jet of the $i$th-order averaging function.  
According to \eqref{1o}, in the case $i=1$, we can write  
\begin{equation}\label{M11}
	\mathcal{M}_1(r)=\sum_{\ell=1}^{N_0}\nu_\ell I_\ell,  
\end{equation}
where $N_0$ is a natural number, $I_\ell = \int_{\Gamma(r)}w_\ell$, with $w_\ell = X_\ell(x,y) dx + Y_\ell(x,y) dy$ being suitable $1$-forms, and $\nu_\ell$ depend linearly on the perturbative parameters of system \eqref{eq:2}, for $\ell = 1,\ldots,N_0$. We note that $w_\ell$ does not depend on these parameters.  
As roughly mentioned above, our strategy in this paper is to find conditions on the perturbative parameters guaranteeing that $M_1^{[j_0]}(r) = 0$ for a certain natural number $j_0$, such that we are able to rewrite \eqref{M11} as  
\begin{equation}\label{M12}
	\mathcal{M}_1(r) = \sum_{\ell=1}^{N_1}\mu_\ell J_\ell,  
\end{equation}  
for some natural $N_1$, such that $J_\ell= \int_{\Gamma(r)}\omega_\ell$ for $\ell = 1,\ldots,N_1$ and each $\omega_\ell$ is in $\Span(w_1,\ldots,w_{N_0})$.  
We are now able to calculate explicitly $J_\ell,$ proving that all vanish. That is, knowing only some coefficients of the Taylor expansion of the first-order averaging function in the neighborhood of the origin up to a certain order and imposing these coefficients equal to $0$, we are able to show that the first-order averaging function is identically zero. 
This approach allows us to analyze cases where we are not able to explicitly calculate \eqref{M11}. 
Therefore, we can proceed to the study of the second-order averaging function.

Our main results are the following ones.  
\begin{theorem}\label{lh1}  
	For the quadratic systems $LV$ and $S_1$, the function $\mathcal{M}_1$ is identically zero if and only if the jet $M_{1}^{[3]}$ is identically zero, and $M_{1}^{[1]} \equiv 0$ does not imply $M_1^{[3]} \equiv 0$.  
\end{theorem}  
\begin{theorem}\label{lh2}  
	For the quadratic systems $S_2$, $S_3$, and $S_4$, the function $\mathcal{M}_1$ is identically zero if and only if the jet $M_{1}^{[5]}$ is identically zero, and $M_{1}^{[3]} \equiv 0$ does not imply $M_1^{[5]} \equiv 0$.  
\end{theorem}  
\begin{theorem}  
	For the quadratic Hamiltonian system $\mathcal{H}$, set  
	\begin{equation}\label{determinant}  
		\begin{aligned}  
			d =& \alpha^3 \beta - \alpha \beta^3 + 6 \alpha^2 \beta \gamma - 2 \beta^3 \gamma + 9 \alpha \beta \gamma^2 + 2 \alpha^3 \delta \\
			&  - 6 \alpha \beta^2 \delta + 9 \alpha^2 \gamma \delta - 9 \beta^2 \gamma \delta - 27 \gamma^3 \delta - 9 \alpha \beta \delta^2 + 27 \gamma \delta^3.  
		\end{aligned}  
	\end{equation}  
	If $d \neq 0$, then $\mathcal{M}_1$ is identically zero if and only if the jet $M_{1}^{[5]}$ is identically zero, and $M_{1}^{[3]} \equiv 0$ does not imply $M_1^{[5]} \equiv 0$.  
	If $d = 0$ and $(\alpha + 3 \gamma)^2 + (\beta + 3 \delta)^2 > 0$, then $\mathcal{M}_1$ is identically zero if and only if the jet $M_{1}^{[3]}$ is identically zero, and $M_{1}^{[1]} \equiv 0$ does not imply $M_1^{[3]} \equiv 0$.  
	If $d = \alpha + 3 \gamma = \beta + 3 \delta = 0$, then $\mathcal{M}_1$ is identically zero if and only if the jet $M_{1}^{[1]}$ is identically zero.  
\end{theorem}  

\begin{theorem}\label{lh4}  
	For the reversible cubic system $CR_1$, the function $\mathcal{M}_1$ is identically zero if and only if the jet $M_{1}^{[7]}$ is identically zero, and $M_1^{[5]} \equiv 0$ does not imply $M_1^{[7]} \equiv 0$.  
\end{theorem}  

The above results are better understood if one recalls that the coefficient of each monomial of even degree $2 j$ in the Taylor expansion of $\mathcal{M}_i$, for any $i$, is contained in the ideal generated by the coefficients of the monomials with degree less than $2 j$, see \cite{CimGas2020, RomSha2009}.  
This explains why in the statements above only even jets are considered.  
In particular, the coefficient of the monomial of degree $2 j$ is zero provided the coefficients of the monomials with less degree are zero.  
Therefore, from the above theorems, we immediately get the following, probably not new, results:  

\begin{corollary}\label{cor1}  
The first-order averaging method applied to the systems $LV$ and $S_1$ provides $1$ limit cycle bifurcating from the origin.  
\end{corollary}  

\begin{corollary}  
The first-order averaging method applied to the systems $S_2$, $S_3$, and $S_4$ provides $2$ limit cycles bifurcating from the origin.  
\end{corollary}  

\begin{corollary}  
For the system $\mathcal{H}$, if $d \neq 0$, then the first-order averaging method provides $2$ limit cycles bifurcating from the origin.  
If $d = 0$ and $(\alpha + 3 \gamma)^2 + (\beta + 3 \delta)^2 > 0$, then $1$ limit cycle bifurcating from the origin is provided.  
If $d = \alpha + 3 \gamma = \beta + 3 \delta = 0$, then no limit cycles bifurcating from the origin are provided.  
\end{corollary}  

\begin{corollary}\label{cor3}  
The first-order averaging method applied to the systems $CR_1$ provides $3$ limit cycles bifurcating from the origin.  
\end{corollary}  

Now, after applying the conditions found in the proofs of Theorems~\ref{lh1} to \ref{lh4}, in order to annihilate $\mathcal{M}_1$, we can proceed with the second-order averaging and analyze the expansion of $\mathcal{M}_2$ to find isolated zeros.  
Some of the following results are again not new.  
Anyway, we (re)obtain them almost directly because, after applying these mentioned conditions guaranteeing that $\mathcal{M}_1 \equiv 0$, we only need to analyze expansions of $\mathcal{M}_2$ to find simple zeros of it.  

\begin{theorem}\label{lh3}  
The second-order averaging method applied to the quadratic systems $LV$, $S_1$, $S_2$, $S_3$, and $S_4$ provides at least $2$ limit cycles bifurcating from the origin.  
	For the system $\mathcal{H}$, at least $2$ limit cycles bifurcate from the origin if $\left(\alpha + 3 \gamma\right)^2 +  \left(\beta + 3 \delta\right)^2 > 0$, and at least $1$ if $\alpha + 3 \gamma = \beta + 3 \delta = 0$.  
\end{theorem}  

\begin{theorem} \label{lh5}
The second-order averaging method applied to the reversible cubic system $CR_1$ provides at least $8$ limit cycles bifurcating from the origin for all but finitely many non-zero $\alpha \in \mathbb{R}$, and at least $6$ if $\alpha = 0$.  
\end{theorem}  
The number of limit cycles obtained here for $\alpha \neq 0$ coincides with the one obtained by \cite{LiZha2014} for $\alpha = 0$.  
But as we will see, the expansion in Taylor series of the second-order averaging function for $\alpha = 0$ does not agree with the formula found in \cite{LiZha2014}.  
In particular, a higher-order analysis is needed to obtain $8$ limit cycles when $\alpha = 0$.  
This is an example where the analysis of the study of the cyclicity by using a family shows that a lower-order expansion of the displacement map provides the same or even better lower bound for it.  
See more details and examples of this fact in \cite{GouTorreGine2021}.

There are other approaches to deal with this problem of rewriting \eqref{M11}. For example, among many others, in \cite{Ili1998, Zol1994, Zol1995} the authors look for $\widetilde{w}_\ell$ being linearly independent in \eqref{M12}.  

We organize the content of this paper as follows: In Section~\ref{section:averaging}, we recall the results for finding limit cycles and explain our method in detail.  
In Section~\ref{section:3}, we prove Theorems~\ref{lh1} and \ref{lh2}.  
Finally, in Section~\ref{section:4}, we prove Theorems~\ref{lh3} and \ref{lh5}.

\section{Limit cycles of polynomial planar systems by perturbing a center}\label{section:averaging}

\subsection{The averaging functions}

We begin by applying the polar change of coordinates $(x,y)=(r\cos\T, r\sin\T)$ to vector field \eqref{eq:2}, obtaining, after taking $\T$ as the new independent variable, the following equivalent differential equation:  
\begin{equation}\label{eq:drdt}
	r'(\T)=\frac{dr}{d\theta}=F_0(\T,r) +  \varepsilon F_1(\T,r) + \varepsilon^2 F_2(\T,r) +\CO(\varepsilon^{3}),
\end{equation}  
with $F_i: [0, 2\pi]\times (0,\rho^*)\rightarrow\R$ analytic functions for small enough $\rho^*$ (because $Z$ is a center), that are $2\pi$-periodic in the variable $\T$, for $i=0,1,2$.  
We denote by $\mathcal{L}(\theta,r,\varepsilon)$ the flow of \eqref{eq:drdt} and, in fact, the $2$-jet of the solution of the initial value problem  
\[
z'(s) = F_0(s, z(s)) +  \varepsilon F_1(s, z(s)) +  \varepsilon^2 F_2(s, z(s)),\quad z(0) = r,  
\]
is given by  
\[
\mathcal{L}_0(\T,r) +  \varepsilon \mathcal{L}_1(\T,r) +  \varepsilon^2 \mathcal{L}_2(\T,r).
\]
We will compute $\mathcal{L}_0,\mathcal{L}_1,$ and $\mathcal{L}_2$ order by order.  
Beginning with $\mathcal{L}_0$, the solution of the initial value problem  
\[
z'(s) = F_0(s,z(s)),\quad z(0) = r,
\]
we obtain that $\mathcal{L}_i$, $i=1,2$, are given by the integral equations  
\begin{equation}\label{y1y2y3}
	\begin{aligned}
		\mathcal{L}_1(\theta,r) &= \int_0^\theta \left(F_1(s,\mathcal{L}_0(s,r)) + \frac{\partial F_0}{\partial r}(s,\mathcal{L}_0(s,r)) \mathcal{L}_1(s,r) \right) ds, \\
		\mathcal{L}_2(\theta,r) & = \frac{1}{2} \int_0^\theta\left(2 F_2(s,\mathcal{L}_0(s,r)) + 2 \frac{\partial F_1}{\partial r}(s, \mathcal{L}_0(s,r)) \mathcal{L}_1(s,r) \right. \\ 
		& \phantom{={}} \left. + \frac{\partial^2 F_0}{\partial r^2}(s,\mathcal{L}_0(s,r)) \mathcal{L}_1(s,r)^2 + \frac{\partial F_0}{\partial r}(s,\mathcal{L}_0(s,r)) \mathcal{L}_2(s,r) \right)ds.
	\end{aligned}
\end{equation}
These expressions, as well as the higher-order ones, are given in \cite{LliNovTei2014}.

We remark that in the following definition, we will use the equivalence between Melnikov and averaging functions, see \cite{Bui2017}.  
So, taking into account \eqref{1o}, up to a multiplicative constant, the \emph{$i$th-order averaging function} is  
\begin{equation*}\label{media}
	\mathcal{L}_i(2\pi,r), 
\end{equation*}
with $i=1,\ldots$  
Evidently, $\mathcal{L}_0(2\pi,r)= 0$, so we do not need $\mathcal{M}_0$ because the unperturbed system is a center. The next proposition gives the $i$th-order averaging method.  
\begin{proposition}\label{hh}
	Assume that $\mathcal{M}_k(r) \equiv 0$, for $k=1,\ldots,i-1$, and $\mathcal{M}_i(r)\not\equiv 0$.  
	If $\mathcal{M}_i$ has $\tau$ simple zeros, then for all $\varepsilon \ne 0$ small enough, the system \eqref{eq:2} has at least $\tau$ limit cycles.  
\end{proposition}

\subsection{Taylor series of the $i$th-order averaging functions}\label{34} 

In order to simplify calculations, we assume further that $Z$ is a non-degenerate center.  
Thus, the functions $ F_i $, for $ i = 0,1,2 $, in \eqref{eq:drdt} are given by  
\begin{equation}\label{F_k}
	F_0(\theta,r) = \frac{r^2 f_0(\theta,r)}{1 + r g_0(\theta,r)}, \quad F_i(\theta,r) = \frac{f_i(\theta,r)}{\left(1+r g_0(\theta,r)\right)^{\eta_i}},
\end{equation}
with $\eta_i\in \mathbb{N}^* $, where $ f_0, g_0 $, and $ f_i $ are polynomials in $ r, \cos\theta, \sin\theta $, for $ i = 1,2 $.  
We expand $ \mathcal{L}_i(\theta, r) $ in $ r $ as  
\begin{equation}\label{varphi} 
	\mathcal{L}_i(\theta,r) = L_{i,j}(\theta,r) + \mathcal{O}\left(r^{j+1}\right), \quad \text{with} \quad L_{i,j}(\theta,r) = \sum_{k=1}^j l_{i,k}(\theta) r^k,
\end{equation}
for $ i = 0,1,2 $, where $ l_{0,1} = 1 $, and each $ l_{i,k} $, $ k = 1,\ldots,j $, corresponds to the coefficient functions of the Taylor expansion in $ r $ of $ \mathcal{L}_i(\theta, r) $.  
The natural number $ j $ denotes the order of the expansion.  
The idea now is to find the functions $ l_{i,k}(\theta) $ iteratively as follows.

We begin with $ \mathcal{L}_0 $.  
By applying the formula of $ F_0 $, given in \eqref{F_k}, to equation \eqref{eq:drdt}, we obtain  
\begin{equation*}
	\Big(1 + \mathcal{L}_0(\theta, r) g_0\left(\theta, \mathcal{L}_0(\theta,  r)\right) \Big){\mathcal{L}_0}'(\theta, r) - \mathcal{L}_0(\theta, r)^2 f_0\left(\theta,\mathcal{L}_0(\theta, r)\right) = 0.
\end{equation*} 
Then we apply the first equation of \eqref{varphi}, and recursively get  
\begin{equation*}
	\begin{aligned}
		{l_{0,2}}'(\theta)& = h_2\left(\cos\theta,\sin\theta\right),\\
		{l_{0,k+1}}'(\theta)& = h_{k+1}\left(\cos\theta,\sin\theta,l_{0,2}(\theta),\ldots, l_{0,k}(\theta)\right),\\ 
	\end{aligned}
\end{equation*}
for $ k=2,\ldots, j-1 $, where the $ h_k $'s are suitable polynomial functions.  
So, after simple integration of trigonometric functions and using the condition $ l_{0,k}(0) = 0 $, we obtain the expressions for the functions $ l_{0,k} $.

To find the Taylor coefficients $ l_{1,k}(\theta) $ of $ \mathcal{L}_1(\theta, r) $, we consider the differential version of the first equation in \eqref{y1y2y3}:
\begin{equation*}
	{\mathcal{L}_1}'(\theta,r) - \left(F_1\left(\theta,\mathcal{L}_0(\theta,r)\right) + \frac{\partial F_0}{\partial r}\left(\theta,\mathcal{L}_0(\theta,r)\right)\mathcal{L}_1(\theta, r) \right) = 0, 
\end{equation*}  
where $ F_0 $ and $ F_1 $ are given in \eqref{F_k}.  
Then, by applying the second equation of \eqref{varphi}, we recursively obtain  
\begin{equation*}\label{orden1}
	\begin{aligned}
		{l_{1,1}}'(\theta) & = h_1\left(\cos\theta,\sin\theta\right),\\ 
		{l_{1,k+1}}'(\theta)&=h_{k+1}\left(\cos\theta,\sin\theta,l_{1,1}(\theta),\ldots,l_{1,k}(\theta)\right), 
	\end{aligned}
\end{equation*}
for $ k=1,\ldots, j-1 $, where again the $ h_k $'s are suitable polynomial functions in $ \cos \theta $, $ \sin \theta $, and the coefficients of the perturbative polynomial. As above, we integrate and obtain the expressions of $ l_{1,k}(\theta) $.

Proceeding in the same way with the second equation of \eqref{y1y2y3} and using \eqref{varphi}, we can obtain the Taylor series coefficients of $ \mathcal{L}_2 $.  
We observe that all these expressions are obtained in a simple manner, since only integration of trigonometric polynomials is necessary.

Now, we define  
\begin{equation}\label{mik}
	m_{i,k} = l_{i,k}(2\pi), 
\end{equation}
for $ i = 1,2 $, and $ k = 1,\ldots, j $, obtaining \eqref{Mism}, the Taylor series of the $i$th-order averaging function.

\section{From local to global}\label{section:3}
The proofs of Theorems~\ref{lh1}, \ref{lh2}, and \ref{lh4} follow directly from the results of this section. We study each of the centers given in \eqref{eq:2_1} and \eqref{eq:2_2} separately in the sequel. 

Assuming the notations of Section~\ref{section:averaging}, our perturbations $Z_i = (P_i, Q_i)$ according to \eqref{eq:2} are by polynomials of the form 
\[
P_i(x,y) =  \sum_{1 \leq k+\ell \leq n} a_{ik\ell} x^k y^\ell, \quad Q_i(x,y)  = \sum_{1 \leq k+\ell \leq n} b_{ik\ell} x^k y^\ell, 
\] 
for $i = 1, 2$, where $n$ is the degree of the original unperturbed center. 

Bellow, we deal with the perturbation of each family separately: each one presents different difficulties, so we need adequate strategies to treat them. 

\begin{proposition}
For the perturbation of system $LV$, the jet $M_{1}^{[3]} \equiv 0$ if and only if 
\begin{equation}\label{cLv}
a_{110} = -b_{101}, \quad a_{102} = - b_{102} - a_{120} - b_{120}.
\end{equation} 
Moreover, this is equivalent to the annihilation of the first-order averaging function \eqref{M11}. 
\end{proposition}
\begin{proof} 
We consider the perturbed system \eqref{eq:2} where $Z$ is the $LV$ system, with $\ell = 1$ and $n = 2$. 
Following the steps \eqref{varphi} to \eqref{mik} of the above algorithm with $j = 3$, we get 
\begin{equation*}
\begin{aligned}
m_{1,1}& =\pi\left(a_{110} + b_{101} \right), \\ 
m_{1,3} & = \frac{\pi}{8}\left(2 a_{120} + 2 a_{102} - a_{110} + 2 b_{120} + 2 b_{102} - b_{101}\right), 
\end{aligned}
\end{equation*}
where $m_{i, j}$ is defined in \eqref{mik}. 
Condition \eqref{cLv} is precisely the solution of $m_{1,1} = m_{1,3} = 0$, which gives $M_{1}^{[3]}=0$. 

We now assume these conditions and, by using formula \eqref{1o} and knowing that 
\[
H(x,y) = x + y - \ln(x+1)(y+1), \quad R(x,y) = (x + 1) (y + 1),
\]
we get expression \eqref{M12} is written as 
\begin{equation*}
\mathcal{M}_1(r) =  b_{110} J_1 + b_{111} J_2 - a_{101} J_3 - a_{111} J_4 + b_{101} J_5 + b_{120} J_6 - a_{120} J_7 + b_{102} J_8,
\end{equation*}		
where $J_i = \int_{\Gamma(r)}\omega_i$ and 
\begin{equation*}
\begin{aligned}
 \omega_1 &= R^{-1} x dx, & \omega_2 &= R^{-1} x y dx, & \omega_3 &= R^{-1} y dy, \\
\omega_4 &= R^{-1} x y dy, &\omega_5 &= R^{-1} \left(y dx + x dy\right), & \omega_6 &=  R^{-1} \left(x^2 dx + y^2 dy\right), \\
\omega_7 &= R^{-1} \left(x^2-y^2\right) dy, & \omega_8 &= R^{-1} y^2 \left(dx+dy\right).
\end{aligned}
\end{equation*}		
Our aim is to prove that $J_i = 0,$ for all $i = 1, \ldots, 8$. 

We first observe that $dH = x(x+1)^{-1} dx + y (y+1)^{-1} dy$. 
We can then write 
\[
R^{-1} x dx = (y+1)^{-1} d H - y (y+1)^{-2} dy, \quad R^{-1} y dy = (x+1)^{-1} d H - x (x+1)^{-2} d x. 
\]
Therefore, since in general, for any $\mathcal{C}^1$ functions $f_1, f_2: A\subset \R^2 \to \R$, we have that $\int_{\gamma} f_1 d f_2 = 0$ if $\gamma$ is a curve contained in a level of $f_2$, and any form $f_1(x) dx$ or $f_1(y) d y$ is exact, it follows that $J_1 =J_2 = J_3 = J_4 = 0$. On the other hand, since $\omega_5, \omega_6,$ $\omega_7 + \omega_8,$ and $H$ are invariant by the change $(x,y)\rightarrow(y,x)$, we obtain $J_5 = J_6 = J_7 + J_8 = 0$ as well. 

Hence, the proof will be finished showing that $J_7 = 0$. In order to proceed, it is convenient to consider in the calculations the first integral $\ee^{-H} = (x+1)(y+1)\ee^{-(x+y)}$ instead of $H$. 
Then we apply to the problem the linear change $(u,v) = (x + y, x - y)$, so that the (new) first integral, after multiplying by $-4$ and adding $4$, is written as 
\[
\widetilde{H} = 4 - (u+2)^2\ee^{-u} + v^2\ee^{-u}, 
\] 
and the $1$-form $\omega_7$ is transformed into
\[
\widetilde{\omega} = 2 u v \left( (u + 2)^2 - v^2\right)^{-2}(du - dv).
\] 
Since $\widetilde{H}(u,-v) = \widetilde{H}(u,v)$ and $\widetilde{\omega}(u,-v) = - \widetilde{\omega}_1(u,v) du + \widetilde{\omega}_2(u,v) dv$, it suffices to prove that, for each $h > 0$ small enough, 
\[
I = \int_{\widetilde{H} = h} \widetilde{\omega}_1(u,v) du =  \int_{u_0}^{u_1} \frac{u \ee^{-u} \sqrt{(u+2)^2 - (4 - h) \ee^u}}{4-h} du, 
\] 
where $u_0 < 0 < u_1$ are the closest to zero ones satisfying $\widetilde{H}(u_0,0) = \widetilde{H}(u_1,0) = h$ (by the expression of $\widetilde{H}$ it is clear that $u_0 = u_0(h)$ and $u_1 = u_1(h)$ exist), is $0$. 

Define $s = s(u)$ in a small neighborhood of $0$ by $s(u) = \sgn(u) \sqrt{4 - (u+2)^2 \ee^{-u}}$. 
This function is continuous and satisfies $s(0)=0$. 
Since $s'(u) = u (2+u)\ee^{-u}/(2 s(u))$ for $u\neq 0$, and $s'(u) \rightarrow 1$ as $u \rightarrow 0$, it follows that $s$ is a $\mathcal{C}^1$ change of variables close to $0$. 
It takes $u_0$ and $u_1$ to $-\sqrt{h}$ and $\sqrt{h}$, respectively. 
Applying this change to $I$ we get 
\[
(4 -h) I = \int_{-\sqrt{h}}^{\sqrt{h}} u \ee^{-u} \sqrt{(u + 2)^2 - \ee^u (4 - h)} \frac{1}{s'(u)}ds = \int_{-\sqrt{h}}^{\sqrt{h}}2 \sqrt{1 - \frac{4-h}{4-s^2}} s ds = 0, 
\]
because the integrand is an odd function. 
\end{proof}

\begin{proposition}\label{explicit}
For the system $\mathcal{H}$, if $d \neq 0$, then $M_{1}^{[5]} \equiv 0$ if and only if 
\[
a_{110} = -b_{101}, \quad a_{111}=-2 b_{102}, \quad b_{111}=-2a_{120}. 
\] 
If $d = 0$, then $M_{1}^{[5]} \equiv 0$ if and only if 
\begin{equation}\label{intermediario}
a_{110} = -b_{101}, \quad (\alpha + 3 \gamma) \left(a_{111} + 2 b_{102}\right) = -(\beta + 3 \delta) \left(b_{111} + 2 a_{120}\right). 
\end{equation}
In this case, if $\left(\alpha + 3 \gamma\right)^2 +  \left(\beta + 3 \delta\right)^2 > 0$, then $M_1^{[3]} \equiv 0$ is equivalent to $M_1^{[5]} \equiv 0$; and if $\alpha + 3 \gamma = \beta + 3 \delta = 0$, then $M_1^{[1]} \equiv 0$ is equivalent to $M_1^{[5]} \equiv 0$. \\
Moreover, the conditions above are equivalent to the annihilation of the first-order averaging function. 
\end{proposition}
\begin{proof} 
In the proof we will assume $\alpha^2 + \beta^2 + \delta^2 + \gamma^2 > 0$, because otherwise we will have the canonical linear center, and in this case if $m_{1,1} = 0$ then $\mathcal{M} = 0$ (this can be seen by the formula of $\mathcal{M}_1$ in \eqref{desh} below: the integrals of all the $1$-forms in this expression are trivially zero because $\Gamma(r)$ is a circle of radius $r$ centered at $(0,0)$). 

As above, we have $m_{1,1} =\pi \big(a_{110} + b_{101} \big)$. 
Assuming from now on $m_{1,1} = 0$, that is 
\[
a_{110} = -b_{101}, 
\] 
it follows that $m_{1,3}$ and $m_{1,5}$ write as 
\begin{equation}\label{msHc}
\begin{aligned}
m_{1,3} &  =\frac{-\pi}{8}\big(\left(\alpha + 3 \gamma \right) \left(a_{111} + 2 b_{102}\right) + \left( \beta + 3 \delta \right) \left(b_{111} + 2 a_{120}\right)\big),\\
m_{1,5}& = \frac{-\pi}{128} \Big( \big(39 \alpha\delta^{2}  + 27 \gamma\delta^{2}  + 30  \alpha \beta \delta+ 30  \beta \gamma \delta+ 5 \alpha^{3} + 15 \alpha^{2} \gamma + 15 \alpha \beta^{2} + 35 \alpha \gamma^{2} \\ 
& \phantom{={}} + 35 \beta^{2} \gamma + 105 \gamma^{3}\big) \left(a_{111} + 2 b_{102}\right) + \big(117 \delta^{3} + 39 \beta\delta^{2}  + 35 \alpha^{2}\delta  + 30 \alpha \gamma \delta \\ 
& \phantom{={}} + 15 \beta^{2} \delta + 15 \gamma^{2} \delta + 15 \alpha^{2} \beta + 30 \alpha \beta \gamma + 5 \beta^{3} + 35 \beta \gamma^{2} \big) \left( b_{111} + 2 a_{120} \right) \Big). 
\end{aligned}
\end{equation}
Here the system is Hamiltonian, with the Hamiltonian function being 
\[
H(x,y) = \frac{x^2 + y^2 + \delta x^3 + \alpha x^2 y + \beta x y^2 + \gamma y^3}{2}. 
\]
(So $R \equiv 1$). 
Then (recall we are assuming $a_{110} = -b_{101}$) we get 
\begin{equation}\label{desh}
\mathcal{M}_1(r) = b_{110} J_1 + b_{120} J_2 - a_{101} J_3 - a_{102} J_4 + b_{101} J_5 + J_6 + J_7
\end{equation} 
where $J_i = \int_{\Gamma(r)}\omega_i$, for $i = 1, \ldots, 7$, and 
\begin{equation*}
\begin{aligned}
&\omega_1  = x dx,\quad  \omega_2 = x^2dx, \quad \omega_3= y dy, \quad \omega_4 = y^2 dy, \quad \omega_5 = y dx + x dy,\\ & \omega_6 =  b_{102} y^2 dx - a_{111} x y dy, \quad \omega_7  = b_{111} x y dx - a_{120} x^2 dy. 
\end{aligned}
\end{equation*} 
Clearly $J_i = 0$ by the exactness of $\omega_i$, for $i = 1,\ldots, 5$. 
Further, $w_6$ is exact, so in particular $J_6 = 0$, if and only if $a_{111} = -2 b_{102}$. 
Analogously, $w_7$ is exact, so in particular $J_7 = 0$, if and only if $b_{111} = - 2 a_{120}$. 
Therefore, assumptions 
\[
a_{111} = -2 b_{102}, \quad b_{111} = - 2 a_{120}, 
\] 
make $\mathcal{M}_1 \equiv 0$. 
And it is clear by \eqref{msHc} that they make $m_{1,3} = 0$ and $m_{1,5} = 0$, respectively. 
	
That this is the only solution of $m_{1,3} = m_{1,5} = 0$ is equivalent to the determinant $5 \pi^2 d/512$, where $d$ is given in \eqref{determinant}, of the matrix of this system (in the variables $a_{111} + 2 b_{102}$ and $b_{111} + 2 a_{120}$) be different from zero. 
Then under condition $d \neq 0$, which is generic on $\alpha, \beta, \gamma, \delta$, we have that $M_1^{[5]} = 0$ implies $\mathcal{M}_1 \equiv 0$. 
Further, yet under condition $d\neq 0$, the linear system on the perturbative parameters giving by $M_1^{[5]} = 0$ has rank $3$, hence we get two limit cycles bifurcating from the original center $\mathcal{H}$ by using first-order averaging theory. 
	
It turns out anyway that when $d = 0$, and so the above assumptions are not the only solution of $M_1^{[5]} \equiv 0$, we, independently, get that the annihilation of the $5$-jet of $\mathcal{M}_1$ provides the annihilation of $\mathcal{M}_1$. 
Actually, as it is natural, in this case we only need to impose $M_1^{[3]} = 0$ or even $M_1^{[1]} = 0$, depending on $\alpha, \beta, \gamma, \delta$ as stated in the proposition, in order to get the annihilation of $\mathcal{M}_1$. 
We shall prove this right below by assuming 
\[
d=0
\] 
from now on and by dividing the analysis into the following cases. 

\subsection*{Case 1: If $\alpha + 3 \gamma \neq 0$ and $\beta + 3 \delta = 0$ (respectively $\alpha + 3 \gamma = 0$ and $\beta + 3 \delta \neq 0$)} 
We have $d$ is a multiple of $\delta (\alpha + 3 \gamma)$ (respectively of $\gamma (\beta + 3 \delta)$), then $d=0$ is equivalent to $\delta = \beta = 0$ (respectively $\gamma = \alpha = 0$). 
Therefore, by \eqref{msHc}, the assumption $m_{1,3} = 0$ is equivalent to $a_{111} + 2 b_{102} = 0$ (respectively $b_{111} + 2 a_{120} = 0$), what immediately yields $m_{1,5} = 0$, as expected in this case. 
Hence $\omega_6$ (respectively $\omega_7$) is exact and so $J_6 = 0$ (respectively $J_7 = 0$). 
Moreover $H(-x,y) = H(x,y)$ (respectively $H(x,-y) = H(x, y)$) and $\omega_7$ under the change $x\to - x$ keeps its expression (respectively $\omega_6$ under the change $y\to - y$ keeps its expression). 
Therefore $J_7 = 0$ (respectively $J_6 = 0$). 

\subsection*{Case 2: If $(\alpha + 3 \gamma) (\beta + 3 \delta) \neq 0$}
Then $m_{1,3} = 0$ is equivalent to 
\begin{equation}\label{456}
a_{111} + 2 b_{102} = -\frac{\beta + 3 \delta}{\alpha + 3 \gamma} \left(b_{111} + 2 a_{120}\right)
\end{equation}
and $m_{1,5} =  10 d \left(b_{111} + 2 a_{120}\right)/(\alpha + 3 \gamma) = 0$, as expected. 
Now we prove that under these assumptions $J_6 + J_7 = 0$. 
	
We consider the rotations 
\begin{equation}\label{rot}
u = x \cos \Theta - y \sin \Theta,\quad v = x \sin \Theta + y \cos \Theta, 
\end{equation}
for each nonzero $\Theta \in (-\pi/2,\pi/2)$, and rewrite $H$ in the new variables $(u,v)$, keeping the same letter $H$ for this, seeking for conditions on $\alpha, \beta, \gamma, \delta$ that guarantees the existence of an angle $\Theta$ such that 
\begin{equation}\label{symmetry}
H(u,-v) = H(u,v). 
\end{equation} 
We will then use this symmetry in order to prove that the transformed one form $w_6 + w_7$ has integral $0$ along the level curves of $H$ close enough to the origin. 
	
Condition \eqref{symmetry} is equivalent to the following linear system in $\beta$, $\gamma$: 
\[
\left(
\begin{array}{cc} 
\!\!(3 \cos^2 \Theta - 1) \sin \Theta & -3 \cos \Theta \sin^2 \Theta \\
-\cos^2 \Theta \sin \Theta & -\cos^3 \Theta
\end{array}
\right) 
\left(\begin{array}{c}
\beta \\
\gamma 
\end{array} 
\right) 
= 
\left(
\begin{array}{c} 
\!\!\cos \Theta \big( \alpha (3 \cos^2 \Theta - 2) + 3 \delta \sin \Theta \cos \Theta \big)\!\!\\
\sin^2 \Theta \big ( \alpha \cos \Theta + \delta \sin \Theta \big)
\end{array}
\right). 
\]
The matrix of the system has determinant $-2 \cos^3 \Theta \sin \Theta$, so we get exactly one solution for each $\Theta$ in our domain, calculated by simply inversion: 
\begin{equation}\label{dfgh}
\beta = \frac{\alpha + 3 \delta r - 3 \alpha r^2 - 3 \delta r^3}{2 r},\quad  \gamma = \frac{- \alpha - 3 \delta r + \alpha r^2 + \delta r^3}{2}, 
\end{equation}
with $r = \tan \Theta$. 
But since $\alpha$, $\beta$, $\gamma$, $\delta$ are given real numbers, we have to see if equations \eqref{dfgh} in $r$ have a real non-zero solution in common. 
We first observe that no common solution can be zero, because this would imply $\alpha = \delta = 0$, contradicting our assumption. 
Now, after multiplying the first equation by $r$, both equations are polynomials in $r$, with leading coefficients multiple of $\delta$. 
Their resultant, with respect to $r$, is a multiple of $\delta^2 d = 0$ and their first subresultant is a multiple of $\delta^2 (\beta + 3 \delta)^2$. 
So if $\delta \neq 0$ it follows the above polynomials have exactly one real root in common. 
Also, if $\delta = 0$, then $\alpha \neq 0$, because otherwise $d = 0$ would imply $\beta \gamma = 0$ and we would contradict our assumptions. 
So in this case, the above polynomials have leading coefficients multiple of $\alpha$ and resultant of the form $\alpha q(\alpha, \beta, \gamma)$, where 
\[
q(\alpha, \beta, \gamma) = \alpha (\alpha + 3 \gamma)^2  - \beta^2 (\alpha + 2 \gamma), 
\] 
and $d$ is a multiple of $\beta q(\alpha, \beta, \gamma)$. 
Since $\beta \neq 0$, we get $q(\alpha, \beta, \gamma) = 0$ and so the resultant is $0$ as well yielding at least a (possible complex) zero $r_0$ in common. 
This zero is a \emph{real} one because plugging it into the second polynomial we get $r_0^2 = (\alpha + 2 \gamma)/\alpha$, which must be positive otherwise $q(\alpha, \beta, \gamma)$ is not $0$. 

We now apply the rotation to the sum $\omega_6 + \omega_7$, obtaining the form $w(u,v)$ that we have to prove has integral $0$ over $\widetilde{\Gamma}(r)$, the curve obtained from $\Gamma(r)$ by the rotation. 
After erasing from $w(u,v)$ the terms with expressions as $f(v) dv$ (which is exact), $g(u,v) du$ and $h(u,v) dv$ such that $g(u,-v) = g(u,v)$ and $h(u,v) = - h(u,-v)$, because $\int_{\widetilde{\Gamma}(r)} g(u,v) du = \int_{\widetilde{\Gamma}(r)} h(u,v) dv = 0$, we obtain the expression
\[
\begin{aligned}
\cos^3\Theta \big(-u v(a_{111} r^3 - 2 a_{120} r^2 + b_{111} r^2 - a_{111} r + 2 b_{102} r - b_{111}) du & \\ 
+ u^2 (b_{102} r^3 - b_{111} r^2 + a_{111} r - a_{120})dv  \big), & 
\end{aligned}
\]
with $r = \tan \Theta$. 
And it is enough to prove this expression has integral $0$ over $\widetilde{\Gamma}(r)$. 
Plugging $a_{111}$ obtained from \eqref{456}, observing that $(\beta+3\delta)/(\alpha + 3 \gamma) = -1/r$, into this expression, we get 
\[
\cos^3 \Theta \big(b_{102} r^3 - b_{111} r^2 - 2 b_{102} r + a_{120} + b_{111}\big) \left(2 u v du + u^2 dv\right), 
\]
which is an exact form, hence it has integral $0$ over $\widetilde{\Gamma}(r)$. 
This finishes the proof of this case. 

\subsection*{Case 3: If $\alpha + 3 \gamma = 0$ and $\beta + 3 \delta = 0$}
Here we will use the same strategy of rotating adequately by means of \eqref{rot} in order to get the symmetry $H(u,-v) = H(u,v)$. 
By using the right above calculations, we see that under the conditions of the present case, the polynomials coming from the equations \eqref{dfgh} are both multiple of 
\[
\delta r^3 - 3 \gamma r^2 - 3 \delta r + \gamma =  \delta (r^2 - 3) r - \gamma (3 r^2 - 1)
\]
which has at least one non-zero real root because $\delta$ and $\gamma$ are not zero simultaneously. 

The existing condition $H(u,-v) - H(u,v) = 0$ provides the polynomial identity
\[
\left(\delta \sin(3 \Theta) - \gamma \cos(3 \Theta)\right) \left(u^2 - \frac{v^2}{3}\right) v = 0, 
\]	
so $\delta \sin(3 \Theta) - \gamma \cos(3 \Theta) = 0$, and hence $H(u,v) = u^2/2 + v^2/2 + A u^3/2 - 3 A u v^2/2$ for $A = \delta \cos(3 \Theta) + \gamma \sin(3 \Theta) \neq 0$. 
The linear change of variables $u \mapsto u/A$, $v \mapsto v/A$ takes $H$ into a multiple of (we keep the same notation) 
\[
H(u,v) = u^2 + v^2 + u^3 - 3 u v^2.
\]
On the other hand (in the original variables), by the exactness of $y^2 dx + 2 x y dy$ and $2 x y d x + x^2 dy$, in order to finish the proof it is enough to prove that $\int_{\Gamma(r)} y^2 dx = \int_{\Gamma(r)} x^2 d y = 0$. 

Applying the rotation and also the right above linear change of variables on these two integrals, after using the exactness of the forms $f(u) du$ and $g(v) dv$ and the identities (coming from the symmetry over the $u$ axis) $\int_{\widetilde{\Gamma}(r)} f(u,v) du = \int_{\widetilde{\Gamma}(r)} g(u,v) dv = 0$ always when $f(u,-v) = f(u,v)$ and $g(u,v) = - g(u,-v)$, where $\widetilde{\Gamma}(r)$ is the transformed curve, it follows that it remains to prove that the forms 
\[
-(2 \sin \Theta) u v du + \sin^3 \Theta (2 u v du + u^2 dv),\quad -(2 \cos \Theta) u v du + \cos^3 \Theta (2 u v du + u^2 dv)
\]
have integral zero over $\widetilde{\Gamma}(r)$. 
And hence, once more by exactness, it is enough to prove that $\int_{\widetilde{\Gamma}(r)} u v du = 0$ where $\widetilde{\Gamma}(r)$ is the curve $u^2(1 + u) + (1 - 3 u) v^2 = r$, for $r>0$ small enough. 
By using the symmetry over the $v$-axis, it is enough to prove that 
\[
\int_{u_0}^{u_1} u \frac{\sqrt{r - u^2 (1 + u)}}{\sqrt{1 - 3 u}} du = 0, 
\]
where $u_0 < 0 < u_1$ are the closest to zero solutions of $u^2 (1 + u) = r$. 
We consider the function $s = s(u) = u \sqrt{1 + u}$ which is smooth invertible close to $0$, fixing $0$, and takes $u_0$ and $u_1$ to -$\sqrt{r}$ and $\sqrt{r}$, respectively. 
Under this change of variables, the integral becomes
\[
\int_{-\sqrt{r}}^{\sqrt{r}} u \frac{\sqrt{r - u^2 (1 + u)}}{\sqrt{1 - 3 u}} \frac{2 \sqrt{1 + u}}{2 + 3 u} ds = \int_{-\sqrt{r}}^{\sqrt{r}} 2 s\frac{\sqrt{r - s^2}}{\sqrt{4 - 27 s^2}} ds = 0, 
\]	
for $r>0$ small enough, as we wanted. 
Clearly in  this case we do not get limit cycles bifurcating from the origin and only $\mathcal{M}_1^{[1]} \equiv 0$ provides $\mathcal{M}_{1} \equiv 0$. 
\end{proof}

\begin{proposition}
For the reversible system $CR_1$, the jet $M_{1}^{[7]}\equiv 0$ if and only if 
\begin{equation}\label{CR}
a_{110} = -b_{101}, \quad a_{130}=b_{121}, \quad  a_{112} = b_{103} =-\alpha(a_{120} + a_{102})+ 2(\alpha^2 + 1)b_{101}-b_{121}. 
\end{equation} 
Moreover, this is equivalent to the annihilation of the first-order averaging function. 
\end{proposition}
\begin{proof}
As in the previous proof, after making $m_{1,1} = 0$, which is equivalent to take $a_{110} = -b_{101}$, we have  
\begin{equation*}
\begin{aligned}
m_{1,3} & = \frac{\pi}{4} \Big(3a_{130} + a_{112} + 4 \alpha a_{120} + 4 \alpha a_{102}+ 3 b_{103} - 8 (\alpha^2 + 1) b_{101}  + b_{121}\Big),\\
m_{1,5} & =\frac{\pi}{4}\Big((3 \alpha^2 + 1) a_{130} + (\alpha^2 + 1) a_{112} + 4  \alpha^3 a_{120} +4  \alpha^3 a_{102}+ (3 \alpha^2 - 1) b_{103}   \\ 
& \phantom{{}=}- 8 \alpha^2 (\alpha^2 + 1) b_{101}+ (\alpha^2 - 1) b_{121}\Big),\\
m_{1,7} & = \frac{\pi}{16}\Big((12 \alpha^4 + 28 \alpha^2 + 3) a_{130} + (4 \alpha^4 + 28 \alpha^2 + 5) a_{112} + 16 \alpha^5 a_{120} + 16 \alpha^5 a_{102} \\ 
& \phantom{{}=} + (12 \alpha^4 - 28 \alpha^2 - 5) b_{103}- 32 \alpha^4 (\alpha^2 + 1) b_{101}  + (4 \alpha^4 - 28 \alpha^2 - 3) b_{121}\Big).  
\end{aligned}
\end{equation*}
It is simple to conclude that the linear system $m_{1,3} = m_{1,5} = m_{1,7} = 0$ in the variables $a_{130}$, $a_{112}$, $b_{103}$ has exactly one solution, as given in \eqref{CR}, proving the first part of the statement. 
In particular, we conclude that at least $3$ limit cycles can be found by using first-order averaging theory. 

To finish the proof, we will show that the assumptions \eqref{CR}, which we set from now on, suffice to make $\mathcal{M}_1 \equiv 0$. 
Here we have that  
\begin{equation*}
H(x,y) = \frac{x^2 + y^2}{1 - 2 x \left(\alpha + x\right)}, \quad R(x,y) = \frac{\left(1 - 2 x \left( \alpha + x\right) \right)^2}{2}, 
\end{equation*}
and so 
\begin{equation*}
\begin{aligned}
\mathcal{M}_1(r) & = b_{130} J_{1} +b_{110}J_2 +b_{120}J_{3}-a_{111}J_{4} -a_{121}J_5 - a_{101}J_6 +b_{102}J_7 +b_{112}J_{8}   \\
& \phantom{{}=}-a_{103}J_{9} +b_{111}J_{10} + b_{121} J_{11} - a_{102}J_{12} -a_{120}J_{13}+b_{101}J_{14}
\end{aligned}
\end{equation*}	
where $J_j = \int_{\Gamma(r)}\omega_j$, for $j = 1, \ldots, 14$, and 
\begin{equation*}
\begin{aligned}
& \omega_{1}=R^{-1} x^3 dx, \quad \omega_2 = R^{-1} x dx, \quad \omega_{3}=R^{-1} x^2dx, \quad \omega_{4}=R^{-1} xydy, \quad \omega_5 = R^{-1} x^2ydy, \\
& \omega_{6}=R^{-1} ydy, \quad \omega_7=R^{-1} y^2dx, \quad \omega_{8} = R^{-1} xy^2 dx, \quad \omega_{9}=R^{-1} y^3dy, \quad   \omega_{10}=R^{-1} xydx, \\
& \omega_{11} =R^{-1} (x^2 - y^2)(ydx - xdy), \quad  \omega_{12}= R^{-1}\left(\alpha y^2 (y dx - x dy) + y^2dy\right), \\
& \omega_{13}=R^{-1}\left(\alpha y^2 (y dx - x dy)+x^2dy \right), \quad \omega_{14}=R^{-1}\left (2 (\alpha^2 + 1) y^2 (y dx-x dy) + ydx + xdy \right). 
\end{aligned}
\end{equation*} 
Clearly $J_1 = J_2 = J_3 = 0$ by exactness. 
In order to deal with the other eleven integrals we apply, in sequence, the changes (i) 
\[
(x,y)= \left(\frac{u}{F + \alpha u}, \frac{v}{F +\alpha u}\right), 
\] 
with $F = \sqrt{1 + (\alpha^2 + 2)u^2}$, yielding 
\[
H = u^2 + v^2, \quad R = \frac{(F + \alpha u)^{-4}}{2},
\]
and (ii) 
\[
(u,v) = (r\cos \theta,r\sin \theta). 
\]
The composition of (i) and (ii) so transforms the curve $\Gamma(r)$ into the parametrized curve $\theta \mapsto \left(r^2,\theta\right)$, $\theta \in (-\pi, \pi)$. 
Further, applying this composition of changes we get the relations 
\[
dx = G_1(r, \theta) dr + \frac{- r \sin \theta}{F (F + \alpha u)^2}  d\theta, \quad d y =  G_2( r, \theta)dr + \frac{u \left(1+ (\alpha^2 + 2) r^2 \right) + \alpha r^2 F}{F (F + \alpha u)^2} d \theta, 
\] 
for suitable functions $G_i: \R^2 \to \R$, $i = 1,2$, and $F = F(u)$, with $u = r \cos \theta$. 

Observe that in the calculation of the transformed integrals, the functions $G_1$ and $G_2$ are irrelevant, because the first coordinate of the transformed curve is constant in $\theta$. 
The integrals $J_4$,$J_5$,...,$J_{9}$ take the form 
\[
\int_{-\pi}^{\pi} K_j(r,\cos \theta) \sin \theta d\theta, 
\] 
for suitable functions $K_j: \R^2 \to \R$, $j=4,\ldots, 9$, and hence they are all zero. 
The integral $J_{10}$ turns to $\int_{-\pi}^{\pi} K_{10}(r,\theta) d \theta$ with $K_{10}(r,\theta) = - r^3 \sin^2 \theta \cos \theta/F$. 
This integral is zero because $K_{10}(r,\theta + \pi) = - K_{10}(r,\theta)$. 
In order to deal with the integrals $J_{11}$, $J_{12}$, $J_{13}$, and $J_{14}$ after applying the composition of transformations, the following identity is useful 
\[
\int_{\Gamma(r)} G(x,y) (y dx - x dy) = \int_{-\pi}^{\pi} G\left(\dfrac{u}{F + \alpha u},\dfrac{v}{F + \alpha u}\right) \frac{-r^2}{(F + \alpha u)^2} d \theta, 
\]
for any function $G: \R^2 \to \R$. 
So the integral $J_{11}$  turns to $\int_{-\pi}^{\pi}-r^4 (\cos^2 \theta - \sin^2 \theta) d \theta = 0$. 
As for $J_{12}$, it is now simple to conclude it turns to $\int_{-\pi}^{\pi} K_{12}(r,\sin^2 \theta) \cos \theta d \theta$, for a suitable function $K_{12}: \R^2 \to \R$, which is zero by the same reason than $J_{10} = 0$. 
Now the integral $J_{13}$ turns to $\int_{-\pi}^{\pi} K_{13}(r,\theta) d\theta$, with 
\[
K_{13}(r,\theta) = \alpha r^4 \cos 2 \theta  + \frac{1+(\alpha^2 + 2) r^2}{F} r^3 \cos^3 \theta. 
\]
Since $K_{13}(r,\theta + \pi) = 2 \alpha r^4 \cos 2 \theta  - K_{13}(r,\theta)$ it readily follows that $J_{13} = 0$. 
Finally, $J_{14}$ turns to $\int_{-\pi}^{\pi} K_{14}(r,\theta) d \theta$, with a similar $K_{14}$, and hence $J_{14} = 0$ as well. 
\end{proof}

For the isochronous systems $S_1$, $S_2$, $S_3,$ and $S_4$ we choose only to make a detailed presentation for $S_4$. 
For $S_1$, $S_2,$ and $S_3$, the proofs are so similar to this that we only deliver a sketch of each of them. 

\begin{proposition}
For system $S_4$, the jet $M_{1}^{[5]} \equiv 0$ if and only if 
\begin{equation}\label{cS4}
a_{110} = -b_{101}, \quad a_{111}=8 b_{101} + 4 b_{120}, \quad b_{102}=\frac{8}{3}b_{101} + \frac{1}{2} b_{120},
\end{equation}
Moreover, this is equivalent to the annihilation of the first-order averaging function. 
\end{proposition}
\begin{proof} 
By using the algorithm in Section~\ref{section:averaging} as above, we get the first coefficients of the $5$-jet of $M_1$ are $m_{1,1} = \pi\left(a_{110} + b_{101} \right)$,  $m_{1,3}  = \pi\left(9 a_{111} - 12 b_{102} - 40 b_{101} - 30b_{120}\right)/9,$ and $ m_{1,5} = 40 \pi\left(21 a_{111} - 48 b_{102} - 40 b_{101} - 60 b_{120}\right)/81.$ 
It is simple to conclude that $m_{1,1} = m_{1,3} = m_{1,5} = 0,$ providing the relations in \eqref{cS4}. 

Now we prove that $\mathcal{M}_1(r)\equiv 0$ under these conditions. 
Here  
\[
H(x,y) = \frac{9 x^2 + (3 + 4 y)^2 y^2} {(3 + 8 y)^2}, \quad R(x,y) = \frac{(3 + 8 y)^5}{54}, 
\]
and so 
\[
\mathcal{M}_1(r) = -a_{101} J_1 - a_{102} J_2 + b_{101} J_3 - a_{120} J_4 + b_{110} J_5 + b_{111} J_6 + b_{120} J_7 
\] 
with $J_i = \int_{\Gamma(r)} \omega_i$, where 
\[
\begin{aligned}
& \omega_1 = R^{-1} y dy, \quad  \omega_2 = R^{-1} y^2 dy, \quad \omega_3 = R^{-1}\left(\frac{(3 + 8 y) y }{3} dx + (1 - 8 y) x dy\right), \quad \omega_4 = R^{-1} x^2 dy, \\ 
& \omega_5 = R^{-1} x dx, \quad \omega_6 = R^{-1} x y dx, \quad \omega_7 = R^{-1} \left(\frac{2 x^2 + y^2}{2} dx - 4 x y dy\right). 
\end{aligned}
\]
By exactness it follows that $J_1 = J_2 = J_3 = 0$, as $d(18 x y/(3 + 8 y)^4) = \omega_3$. 
After applying the (canonical) change 
\[
(u,v) = \left(3 x S^2, (3 + 4 y) y S^2\right), 
\] 
with $S = (3 + 8 y)^{-1}$, 
it follows that $H$ is transformed into $u^2 + v^2$. 
From the relations 
\[
R = \frac{1}{54 S^5}, \quad x = \frac{u}{3 S^2}, \quad y = \frac{1 - 3 S}{8 S}, \quad S = \frac{\sqrt{1 - 16 v}}{3}, 
\] 
by further considering the polar change $(u, v) =  (r\cos\T, r\sin\T)$, 
we get 
\[
d x = G_1(r, \theta) dr +\left( \frac{16 u^2}{27 S^4} - \frac{v}{3 S^2}  \right) d \theta, \quad dy = G_2(r, \theta) d r + \frac{u}{9 S^3} d\theta
\]
with suitable functions $G_1, G_2: \R^2 \to \R$. 
Then, by also considering that the transformed curve can be parametrized by $(r^2, \theta)$, $\theta \in (-\pi, \pi)$, it is simple to see by inspection that the transformed integrals $J_4$, $J_5$ and $J_6$ are of the form 
\[
\int_{-\pi}^{\pi} K_i(r,\sin \theta) \cos \theta, 
\]
with suitable functions $K_i: \R^2 \to \R$, $i = 4,5,6$. 
But integrals of the form $\int_{-\pi}^{\pi}g(\theta) d\theta$, where $g$ is a $2 \pi$-periodic function satisfying $g(\pi/2 + \theta) = - g(\pi/2 - \theta)$, are zero, because $\int_{-\pi}^{\pi} g(\theta) d\theta = \int_{-\pi}^{\pi} g(\pi/2 + \theta) d\theta$ by the periodicity, and this is zero as $g(\pi/2 + \theta)$ is odd. 
Therefore, $J_4 = J_5 = J_6 = 0$. 
Finally, we turn our attention to $J_7$. 
After the transformations, this integral is written $\int_{-\pi}^{\pi} K_7(r,\theta)$, with 
\begin{equation*}
\begin{aligned}
K_7(r,\theta) & = \frac{r}{15552 S^3}\Big(16 r \cos^2 \theta \left(9 S^2 (3 S -1) + 4 r (16 r - \sin \theta)\right) \\ 
& \phantom{={}} +9 S^2 \sin \theta \left(9 S^2 (3 S - 1) - 8 r (16 r - \sin \theta )\right) \Big). 
\end{aligned}
\end{equation*}	
It is then a matter of some calculations to see that the following expression is a primitive of $\theta \mapsto K_(r, \theta)$: 
\begin{equation*}
\frac{r \cos \theta}{1728 R} \Big( 8 r (16 r - \sin \theta) - 9 S^2 (3 S - 1)\Big). 
\end{equation*}	
Since this is a $2 \pi$-periodic function, it follows that $J_7 = 0$. 
\end{proof}

Now it follows the sketch for $S_1$, $S_2$, and $S_3$. 

\subsection*{The case $S_1$} 
The jet $M_1^{[3]}\equiv 0$ if and only if 
\[
a_{110} = -b_{101} = -\frac{b_{102} + b_{120}}{2}. 
\]
It follows that $\mathcal{M}_1$ is zero if and only if these conditions hold, because here 
\[
H(x,y) = \frac{x^2 + y^2}{1 + 2 y}, \quad R(x,y) = \frac{(1 + 2 y)^2}{2}, 
\]
and so
\[ 
\mathcal{M}_1(r) = -a_{101} J_1 - a_{102} J_2 + b_{102} J_3 - a_{120} J_4 + b_{110} J_5 + b_{111} J_6 - a_{111} J_7 + b_{120} J_8 , 
\]
where $J_i = \int_{\Gamma(r)} \omega_i$ and 
\[
\begin{aligned}
& \omega_1 = R^{-1} y dy \quad \omega_2 = R^{-1} y^2 d y, \quad \omega_3 = \frac{R^{-1}}{2} \left((1+2 y)y dx + x dy\right), \quad \omega_4 = R^{-1} x^2 dy, \\
&  \omega_5 = L x dx,\quad  \omega_6 = R^{-1} x y dx, \quad  \omega_7 = R^{-1} x y dy, \quad \omega_8 = \frac{R^{-1}}{2}\left((2 x^2 + y) dx + x dy\right),
\end{aligned} 
\]
The forms $\omega_1$, $\omega_2,$ and $\omega_3$ are exact as $d(2 x y/(1+2 y)) = \omega_3$, then $J_1 = J_2 = J_3 = 0$. 
In order to analyze the other integrals, we apply the following change of variables (it is different from the canonical one):
\[
x = u \left(v + \sqrt{1+v^2}\right), \quad y = v \left(v + \sqrt{1 + v^2}\right), 
\]
that carries $H$ into $u^2 + v^2$, and right after we apply the change $(u, v) = (r \cos \theta, r \sin \theta)$, so the transformed curve $\Gamma(r)$ can be now parametrized by $\theta \mapsto (r, \theta)$, $\theta \in (-\pi, \pi)$. 
Then it is a question of a simple inspection, after calculating $dx$ and $dy$ in the new variables, to see that $J_4$, $J_5,$ and $J_6$ get transformed into $\int_{-\pi}^{\pi} K_i(r,\sin \theta) \cos \theta d \theta$, for suitable functions $K_i$, $i=4,5,6$, respectively, hence $J_4 = J_5 = J_6 = 0$, as in the previous proof. 
Further, $J_7$ get transformed into $\int_{-\pi}^{\pi} K_7(r, \cos \theta) \sin \theta d\theta$ which is zero as well and $J_8$ get transformed into $\int_{-\pi}^{\pi} K_8(r,\theta) d \theta$, with 
\[
K_8(r, \theta) = (1 + 2 r^2) (u^2 - v^2)-\frac{2 (1 + r^2) (u^2 - v^2)}{\sqrt{1 + v^2}} v , 
\]
with $u = r \cos \theta$ and $v = r \sin \theta$, and hence $J_8 = 0$. 

\subsection*{The case $S_2$}
The jet $M_1^{[5]}\equiv 0$ if and only if 
\[
a_{110} = -b_{101},\quad  b_{102} = a_{111}, \quad  b_{120} = 0. 
\]
These conditions are equivalent to $\mathcal{M}_1 \equiv 0$ because here 
\[
H = \frac{x^2 + y^2}{(1 + y)^2}, \quad R = \frac{(1+y)^3}{2}, 
\]
and 
\[
\mathcal{M}_1 = - a_{101} J_1 - a_{102} J_2 + a_{111} J_3 + - a_{120} J_4 + b_{110} J_5 + b_{111}J_6 + b_{101} J_7, 
\]
where $J_i = \int_{\Gamma(r)} \omega_i$ and
\[
\begin{aligned}
& \omega_1 = R^{-1} y dy, \quad \omega_2 = R^{-1} y^2 dy, \quad \omega_3 = R^{-1} y\left( y dx - x dy\right), \quad \omega_4 =  R^{-1} x^2 dy, \quad \omega_5 = R^{-1} x dx, \\ 
& \omega_6 = R^{-1} x y dx, \quad \omega_7 = R^{-1} (y dx + x dy). 
\end{aligned} 
\]
$J_1 = J_2 = 0$ because $\omega_1$ and $\omega_2$ are exact. 
To analyze the other integrals, we apply the canonical change 
\[
x = \frac{u}{1-v}, \quad  y = \frac{v}{1-v}, 
\]
carrying $H$ into $u^2 + v^2$, and after the change $(u,v) = (r \cos \theta, r \sin \theta)$, so $\Gamma(r)$ get transformed into the parametrized curve $(r, \theta)$, $\theta \in (-\pi, \pi)$. 
Then $J_3$, $J_4$, $J_5$, $J_6$, and $J_7$ are transformed to $\int_{-\pi}^{\pi} K_3(r, \cos \theta) \sin \theta d \theta$, $\int_{-\pi}^{\pi} K_i(r, \sin \theta) \cos \theta d \theta$, $i = 4,5,6$, and $\int_{-\pi}^{\pi} \frac{\partial K_7}{\partial \theta}(r, \theta)$, for suitable functions $K_3, \ldots, K_6$, and $K_7(r, \theta) = 2 \cos \theta \left( \sin \theta  - r \right) r^2$. 
Hence $J_3 = J_4 = J_5 = J_6 = J_7 = 0$. 

\subsection*{The case $S_3$}
The jet $M_1^{[5]}\equiv 0$ if and only if 
\[
a_{110} = -b_{101} =  \frac{3 b_{102}+ 4 b_{120}}{16}, \quad a_{111} = - \frac{b_{102}}{2}. 
\] 
And these conditions are equivalent to $\mathcal{M}_1 \equiv 0$ because here 
\[
H = \frac{16 x^4 - 24 x^2 y + 9 x^2 + 9 y^2}{3 - 16 y}, \quad R = \frac{(16 y-3)^2}{6(32 x^2 - 24 y + 9)}, 
\]
so 
\[ 
\mathcal{M}_1 = - a_{101} J_1 - a_{102} J_2 - a_{120} J_3 + b_{111} J_4 + b_{110} J_5 + b_{102} J_6 + b_{120} J_7,
\]
where $J_i = \int_{\Gamma(r)} \omega_i$ and 
\[
\begin{aligned}
& \omega_1 = R^{-1} y dy, \quad \omega_2 = R^{-1} y^2 dy, \quad \omega_3 = R^{-1} x^2 dy, \quad \omega_4 = R^{-1} x y dx, \quad \omega_5 = R^{-1} x dx, \\ 
& \omega_6 = \frac{R^{-1}}{16}\big((16 y - 3) y dx + (8 y - 3)x dy\big),  \quad \omega_7 = \frac{R^{-1}}{4}\big((4 x^2 - y)dx - x dy\big). 
\end{aligned}
\]
We apply the canonical change 
\[
x = \frac{3u}{8 v + 1}, \quad y = 3 \frac{4 u^2 + 8 v^2 + v}{(8 v + 1)^2}, 
\]
that transforms $H$ into $g(u^2 + v^2)$ for a suitable $g$, and the further change $(u, v) = (r \cos \theta, r \sin \theta)$, so the curve $\Gamma(r)$ turns to $\theta \mapsto (r, \theta)$, $\theta \in (-\pi, \pi)$. 
The first five integrals above turn to $\int_{-\pi}^{\pi} K_i(r,\sin \theta)\cos \theta d \theta$, for $i = 1,\ldots, 5$, so all are zero. 
And the last two turn to $\int_{-\pi}^{\pi} \frac{\partial K_i(r,\theta)}{\partial \theta} d\theta$, for $i = 6,7$, where 
\[
K_6(r, \theta) = \frac{81 r^2 \cos \theta \left( 4 r (1 + \sin^2 \theta) + \sin \theta \right)}{8 (64 r^2 - 1) (1 + 8 r \sin^2 \theta)^2}, \quad K_7(r, \theta) = - \frac{27 r^2 \cos \theta (8 r + \sin \theta)}{2(64 r^2 - 1)^2} , 
\]
which are $2 \pi$-periodic, and hence $J_6=  J_7 = 0$ as well. 

\smallskip

The proofs of Corollaries~\ref{cor1} to \ref{cor3} follow directly from the independence of the coefficients $m_{1,1}$ and $m_{1,3}$; or $m_{1,1}$, $m_{1,3}$, and $m_{1,5}$; or $m_{1,1}$, $m_{1,3}$, $m_{1,5}$, and $m_{1,7}$, depending on the case analyzed above. 
This independence makes one be able to find one, two or three simple zeros of $\mathcal{M}_1(r)$, respectively. 
Then Proposition~\ref{hh} guarantees the respective systems have at least $1$, $2$, or $3$ limit cycles. 

\section{Applying second-order averaging theory: the proof of Theorem~\ref{lh3}}\label{section:4}
By using the results of Section~\ref{section:3}, we prove here Theorem~\ref{lh3}. 
Precisely, for each system considered in these paper, we set the conditions, founded in Section~\ref{section:3}, annihilating $\mathcal{M}_1$. 
And we consider perturbations of degree $2$ for each of them, considering $\ell=2,$ in $\eqref{eq:2}.$ 
Then we follow the algorithm of Section~\ref{34}, obtaining the first terms of the Taylor expansion of $\mathcal{M}_2(r)$, so that we can apply Proposition~\ref{hh}. 

As said in Section~\ref{se:1}, most of the results below are not new, so we deliver only a sketch of the proofs. 

\subsection{The case $LV$} 
We assume conditions \eqref{cLv} and calculate the $7$-jet 
\[
M_2^{[7]}(r) = m_{2,1} r + m_{2,3}r^3 + m_{2,5} r^5 + m_{2,7} r^7
\] 
of $\mathcal{M}_2$ by following the algorithm of Section~\ref{section:averaging}. 
Then we introduce the parameters $A_1$ and $A_2$ by solving $m_{2,1} = A_1$ and $m_{2,3} = A_2$, in the parameters $a_{210}$ and $b_{202}$, respectively. 
We further introduce new auxiliary parameters $A_3$ and $A_4$ by substituting $b_{102} = A_3 - b_{120} + b_{101}$ and $b_{110} = -2 \left(A_4 - b_{111} - a_{111} - a_{120} - b_{102}\right) - a_{101}$. 
Then  
\begin{equation*}
\begin{aligned}
& m_{2,1} = A_1, \quad m_{2,3} = A_2, \quad m_{2,5} = \frac{11}{36} A_2 - \frac{67}{1280} A_1 - \frac{\pi}{6}  A_3 A_4, \\
&m_{2,7} = \frac{979}{5760} A_2 - \frac{64037}{4354560} A_1 - \frac{53\pi}{288}  A_3 A_4. 
\end{aligned}
\end{equation*}
Since $A_k$, $k =1,2,3,4$, are free parameters, we certainly get at least $2$ simple zeros of $M_2^{[7]}(r)$ and so of $\mathcal{M}_2(r)$. 
Therefore, from Proposition~\ref{hh}, it follows that there exists a quadratic perturbation of second order of system $LV$ having at least $2$ limit cycles. 

Here we went further, calculating up to the $11$-jet of $M_2$, but did not find new independent coefficients to guarantee more isolated zeros. 
Actually, we believe they do not exist. 
	
\subsection{The case $\mathcal{H}$} 
We assume $a_{110} = - b_{101}$, which is a necessary condition for $\mathcal{M}_1$ to be zero, independently on the cases of Proposition~\ref{explicit}. 
Further, it is simple to see that $m_{2,1} = A_1$ if and only if $a_{210} = \left(A_1 - b_{201}\right)/(2 \pi)$, a condition we also assume from now on. 
Here $A_1$ is a new parameter we just introduce. 

Then, independently on the other conditions appearing in Proposition~\ref{explicit}, it happens that for given $A_2, A_3\in \R$, the system of equations $m_{2,3} = A_2$, $m_{2,5} = A_3$ is linear in the variables $(a_{211} + 2 b_{202})$ and $(b_{211} + 2 a_{220})$. 
Actually, this system has the form $2 L X = b$, where $b = (b_1,b_2)^t$ is a column vector with entries $b_i = A_{i+1} +$ suitable homogeneous polynomials in the parameters $a_{1ij}, b_{1ij}$, and $L$ is precisely the same matrix as the one of the linear system $m_{1,3} = m_{1,5} = 0$, in the variables $(a_{111} + 2 b_{102})$ and $(b_{111} + 2 a_{120})$, already appearing in the proof of Proposition~\ref{explicit}, see \eqref{msHc} above. 

So in case $d \neq 0$ ($d$ is given in \eqref{determinant} and is a multiple of the determinant of the matrix $L$), we solve this system $2L X = b$ and we get $a_{211} = - 2 b_{202} + q_1$ and $b_{211} = -2 a_{220} + q_2$, where $q_1$ and $q_2$ are suitable polynomials of degree $2$ in the parameters $a_{1ij}$,$b_{1ij}$, $A_3$, $A_5$. 
Under this assumption on $d$, it follows by Proposition~\ref{explicit} that $a_{111} = - 2 b_{102}$ and $b_{111} = -2 a_{120}$ as well in order to get $\mathcal{M}_1 = 0$. 
We get so far 
\[
\begin{aligned}
& m_{2,1} = A_1, \quad m_{2,3} = A_2, \quad m_{2,5} = A_3, \\
& m_{2,7} = \frac{3}{64} \left(21 \alpha^2 + 14 \alpha \gamma + 21 \beta^2 + 14 \beta \delta + 93 \delta^2 + 77 \gamma^2\right)A_3.
\end{aligned}
\]
We calculate up to $m_{2,9}$ and could not see new free parameters popping up. 
So we confirm at least two limit cycles in case $d\neq 0$. 

On the other hand, when $d = 0$, the above considered system cannot be solved by inverting the matrix $L$.  
We divide the analysis in the two cases given by Proposition~\ref{explicit}, that is, $\left(\alpha + 3 \gamma\right)^2 +  \left(\beta + 3 \delta\right)^2 > 0$, or $\alpha + 3 \gamma = \beta + 3 \delta = 0$. 

For the first one, we first analyze when $\alpha + 3 \gamma \neq 0$. 
Here we can isolate $a_{111}$ in equation \eqref{intermediario} and $a_{211}$ in equation $m_{2,3} = A_2$. 
This will give 
\[
\begin{aligned}
m_{2,5} & = \frac{5\pi \left(-2\alpha^3 - 9 \alpha^2\gamma + 6\alpha\beta^2 + 18\alpha\beta\delta + 9\beta^2\gamma - 81\delta^2\gamma + 27\gamma^3\right)\left(b_{111} + 2 a_{120}\right)}{48(3 \gamma + \alpha)} b_{120} \\
& \phantom{{}=} - \frac{5\pi \left(\alpha^3 + 6\alpha^2 \gamma - 3\alpha\beta^2 - 12\alpha\beta\delta - 9\alpha\delta^2 + 9\alpha\gamma^2 - 6\beta^2\gamma - 18\beta\delta\gamma\right)\left(b_{111} + 2 a_{120}\right)}{16(3\gamma + \alpha)} b_{102}   + q, 
\end{aligned}
\] 
where $q$ depends only on parameters different from $b_{120}$ and $b_{102}$. 
Then if either the coefficient of $b_{120}$ or of $b_{102}$ in the expression of $m_{2,5}$ is not zero, we can solve $m_{2,5} = A_3$ in $b_{120}$ or in $b_{102}$, respectively. 
It turns out that in both cases, after imposing $A_1 = A_2 = A_3 = 0$, $m_{2,7}$ and $m_{2,9}$ become multiples of $d$, which is $0$. 
And we thing this happens for higher coefficients, so we can assure no more than $2$ limit cycles. 
On the other hand, if both the mentioned coefficients above are zero, that happens either if $b_{111} + 2 a_{120} = 0$ or if the polynomials of degree $3$ in $\alpha$ are zero, it follows that: 
(i) if $b_{111} + 2 a_{120} = 0$, it is simple to see that under $A_1 = A_2 = 0$, it happens that $m_{2,5} = m_{2,7} = m_{2,9} = 0$, that is, we only guarantee $1$ limit cycle; 
(ii) if this is not zero but the two polynomials of degree $3$ in $\alpha$ are zero, then the resultant between them, 
\[
108\gamma \left(3\delta + \beta\right)^6 \left(-3\delta^2 + \gamma^2\right),
\] 
must be zero as well. 
If $\delta = 0$, since then $\alpha \neq 0$, it is not difficult to conclude we will arrive to $d \neq 0$. 
If $3\delta+\beta=0$, it follows that we must have $\alpha+3\gamma = 0$, a new contradiction. 
If $\gamma^2 = 3 \delta^2$, that is $\gamma = \pm \sqrt{3} \delta$, then one sees that it is necessary that $\alpha = \pm \sqrt{3} \beta$ and $\delta = \beta$. 
Under these assumptions and also putting $A_1 = A_2 = 0$, it follows that $m_{2,5}$, $m_{2,7},$ and $m_{2,9}$ are multiples of a same expression (with some free parameters) that can be done different from zero. 
So we get at least $2$ limit cycles. 
Since $b_{111}$ and $a_{120}$ are free parameters, we can always assume $b_{111} + 2 a_{120} \neq 0$, and so $2$ limit cycles are unfolded here. 

Now, when $\beta + 3 \delta \neq 0$ a completely analogous situation occurs, but now with $a_{120}$ and $a_{102}$ in place of $b_{120}$ and $b_{102}$, respectively, and $a_{111} + 2 b_{102}$ in place of $b_{111} + 2 a_{120}$. 
So we can also unfold at least two limit cycles here. 

Finally, we analyze the second case of Proposition~\ref{explicit}, that is $\alpha = - 3\gamma$ and $\beta=-3 \delta$. 
Here we have
\[
m_{2,3} = \frac{\pi (a_{111} + 2 b_{102})}{2} a_{102} -\frac{\pi(b_{111} + 2 a_{120})}{2}b_{120} + q, 
\]
where $q$ depends on the other parameters different from $a_{102}$ and $b_{120}$. 
Acting as above, if either $a_{111} + 2 b_{102} \neq 0$ or $b_{111} + 2 a_{120}\neq 0$, we can solve $m_{2,3} = A_2$ in $a_{102}$ or $b_{120}$, respectively. 
In both cases it follows that $m_{2,5}$, $m_{2,7}$ and $m_{2,9}$ are zero when $A_1 = A_2 = 0$. 
So we can guarantee at least $1$ limit cycle. 
Now if $a_{111} + 2 b_{102} = b_{111} + 2 a_{120} = 0$, it follows that $m_{2,3} = m_{2,5} = m_{2,7} = m_{2,9} = 0$ if $A_1 =0$, and hence no limit cycle can be unfolded here. 
Of course, since $a_{111}$, $b_{102}$, $b_{111},$ and $a_{120}$ are free parameter, we assure at least $1$ limit cycle. 

The proof of Theorem~\ref{lh3} is completed by observing that $\alpha + 3 \gamma = \beta + 3 \delta = 0$ imply $d = 0$ as well. 

\subsection{The case $CR_1$}
We assume conditions \eqref{CR} and calculate $M_2^{[17]}$. 
Then we introduce the parameters $A_1, A_2, A_3,$ and $A_4$ and solve the system $m_{2,2k+1} = \pi A_{k+1}$ for $k=0,1,2,3$ in the coefficients $b_{201},$ $b_{221},$ $b_{203},$ and $a_{230}$, obtaining 
\begin{equation}\label{4566}
\begin{aligned}
m_{2,1} & = \pi A_1, \quad m_{2,3} = \pi A_2, \quad m_{2,5} = \pi A_3, \quad m_{2,7} = \pi A_4,\\
m_{2,9} & =\frac{\pi}{4}\Big(3 \alpha^2(96\alpha^4 + 16\alpha^2 + 1)A_2- (352\alpha^4+56\alpha^2+3)A_3+ 4 (17\alpha^2+2)A_4\Big) + Q_9,\\ 
m_{2,11} & = \frac{\pi}{4}\Big( \alpha^2(2592\alpha^6+1488\alpha^4 + 177\alpha^2+7)A_2 -(3104\alpha^6 + 1736\alpha^4+191\alpha^2 + 7)A_3 \\
& \phantom{{}=}+ 2(258\alpha^4+124\alpha^2+7)A_4\Big) + Q_{11},\\
m_{2,13} & = \frac{\pi}{8}\Big(\alpha^2(32832\alpha^8+33696\alpha^6 + 10338\alpha^4 + 974\alpha^2 +27)A_2 - (38976\alpha^8 + 39312\alpha^6\\
& \phantom{{}=}  + 11710\alpha^4 + 1022\alpha^2 + 27)A_3 + 4 (1538\alpha^6 + 1404\alpha^4 + 343\alpha^2 + 12)A_4\Big) + Q_{13},\\
m_{2,15} &=\frac{\pi}{16}\Big( \alpha^2(360576\alpha^{10} + 540480\alpha^8 + 292164\alpha^6 + 62428\alpha^4 + 4824\alpha^2 + 99)A_2\\
&\phantom{{}=} -(426112\alpha^{10}+630560\alpha^8 + 335228\alpha^6 + 69244\alpha^4 + 4989\alpha^2 + 99)A_3\\
&\phantom{{}=}+(65552\alpha^8 + 90080\alpha^6 + 43064\alpha^4 + 6816\alpha^2+165)A_4\Big) + Q_{15},\\			
m_{2,17} & = \frac{\pi}{64}\Big( \alpha^2(7340544\alpha^{12} + 14548224\alpha^{10} + 11576592\alpha^8 + 4278128\alpha^6 + 701856\alpha^4\\
&\phantom{{}=} + 45276\alpha^2 + 715)A_2- (8651264\alpha^{12}+16972928\alpha^{10} + 13354480\alpha^8 + 4842992\alpha^6 \\
&\phantom{{}=}+ 765876\alpha^4 + 46420\alpha^2 + 715)A_3 + 4(327696\alpha^{10} + 606176\alpha^8 + 444472\alpha^6  \\
&\phantom{{}=}+ 141216\alpha^4+ 16005\alpha^2 + 286)A_4\Big) + Q_{17},\\			
\end{aligned}
\end{equation}
where $Q_k$ are suitable homogeneous polynomials of degree two in the variables $a_{1ij}$, $b_{1ij}$, whose coefficients only depend polynomially on $\alpha$. 
Their expressions are too big to be presented here. 
Since $A_k$, $k=1,\ldots,4$, are free parameters, it follows, as above, that so far we can unfold at least $3$ limit cycles from a quadratic perturbation of $CR_1$, independently on $\alpha$. 

But we want to analyze the possible existence of more limit cycles. 
So we first consider $\alpha \neq 0$. 
The strategy is to seek for a point in the space of parameters contained in a \emph{transversal} intersection of the algebraic varieties 
\begin{equation*}
m_{2,1} = m_{2,3} = \cdots = m_{2,15} = 0, 
\end{equation*}
such that $m_{2,17}\ne0$. 
Then it is not difficult to conclude that by slightly moving this point, one can get freedom enough in the coefficients of $M^{[17]}(r)$ in order to produce $8$ simple zeros of $M^{[17]}(r)$, and so of $\mathcal{M}(r)$, getting at least $8$ limit cycles bifurcating from the origin. 
See also \cite[Theorem 3.1]{Chr2005}. 
We now push forward this strategy, and in the process the meaning of ``transversal'' will be clear. 
First make $A_1 = A_2 = A_3 = A_4 = 0$. 
This is the only solution of $m_{2,1} = m_{2,3} = m_{2,5} = m_{2,7} = 0$, and it is ``transversal'' because the map $\R^4\ni(A_1, A_2, A_3, A_4) \mapsto \left(m_{2,1},m_{2,3},m_{2,5},m_{2,7}\right) \in \R^4$ is a diffeomorphism fixing $0$. 
Before continuing, we eliminate some variables by putting $a_{103} = b_{120} = a_{120} = b_{112} = a_{101} = b_{101} = b_{110} = 0$, $a_{121}  =  b_{111} = \alpha$ obtaining the system of $4$ equations 
\[ 
m_{2,2k+1}, \quad k = 4,5,6,7, 
\]
in the $4$ variables $a_{102}, a_{111}, b_{102}, b_{130}$
That is, we have a map 
\[
F: \R^4 \to \R^4, 
\] 
fixing $0$, whose coordinate functions are $m_{2,2k+1}$, $k = 4,5,6,7$, in the variables $a_{102}, a_{111},$ $b_{102}, b_{130}$ and we are looking for the fiber $0$ of $F$. 
By using any algebraic manipulator it is simple to find three solutions. 
One of them is certainly not transversal, and between the other two we choose the following: 
\begin{equation*}
\begin{aligned}	
a_{102}& = \frac{6\alpha^2 + 5 - T}{2\alpha}, \\
a_{111} &= \frac{\alpha^2}{S} \Big( \alpha^2 T \left(48\alpha^{10}-112\alpha^{8} + 236\alpha^6 - 421\alpha^4 - 110\alpha^2 + 30\right) \\
& \phantom{={}} + 704\alpha^{14} - 576\alpha^{12} + 968\alpha^{10} - 2218\alpha^8 - 3255\alpha^6 - 350\alpha^4 + 475\alpha^2 - 700 \Big),\\  
b_{102} &= \frac{\alpha^2}{S} \Big( \alpha^2 T  \left(24\alpha^{10}-56\alpha^{8}+118\alpha^6 - 228\alpha^4 - 55\alpha^2 + 15 \right) \\
& \phantom{={}} + 368\alpha^{14} - 296\alpha^{12} + 492\alpha^{10} - 1014\alpha^8 - 1722\alpha^6 - 105\alpha^4 + 125\alpha^2 - 175 \Big) , \\
b_{130} &= \frac{\alpha^3}{S}\Big(-2 T \left( 16 \alpha^{12} - 32\alpha^{10} + 64\alpha^8 - 128\alpha^6 - 94 \alpha^4 + 20 \alpha^2 - 5\right) \\
&  \phantom{={}} - 480\alpha^{14} + 224\alpha^{12} - 448\alpha^{10} + 896\alpha^8 + 2828\alpha^6 + 574\alpha^4 - 140\alpha^2 + 125\Big), 
\end{aligned}
\end{equation*}
where 
\[
T = \sqrt{(6 \alpha^2 + 5)^2 + 24 \alpha^2}, \text{ and } S = 32\alpha^{16} - 420\alpha^8 + 14\alpha^6 - 140\alpha^4 + 225\alpha^2 - 350. 
\]
Moreover, this solution makes 
\[
m_{2,17} = K \Big(T \left(8 \alpha^8 + 4 \alpha^6 - 2 \alpha^2 + 5\right) - 48 \alpha^{10} - 80 \alpha^8 - 4 \alpha^6 - 12 \alpha^4 + 20 \alpha^2 - 95 \Big) 
\]
where
\[
K = \frac{\alpha^9(\alpha^2 + 2)^2 \pi}{144 \left(32 \alpha^{14} - 64 \alpha^{12} + 128 \alpha^{10} - 256 \alpha^8 + 92 \alpha^6 - 170 \alpha^4 + 200 \alpha^2 - 175\right)}. 
\]
So $m_{2,17}$ is non-zero for all but finitely many $\alpha \in \R$. 
And now comes the transversality of this solution. 
The Jacobian determinant of $F$ at this point is 
\[
\frac{\pi^4 \left(6\alpha^2 + 5 - T\right)^2 \left(\alpha^2 + 2\right)^7 \alpha^8 T}{464486400}, 
\] 
which is a nonzero quantity for all but finitely many $\alpha$. 
Then $F$ is a local diffeomorphism justifying the idea of freedom on the coefficients $m_{2,k}$ mentioned above. 

Here we have calculated up to the $21$-jet and did not find more than $8$ limit cycles. 
		
On the other hand, when $\alpha =0 $, after putting $A_1 = A_2 + A_3 = A_4 = 0$ and rewriting 
\begin{equation*}
	\begin{aligned}
b_{120} &= A_5 - 6 a_{111} + 11 b_{102}, & a_{120} &= A_6 - a_{102}, \\
b_{102} &= \frac{A_7}{12} - \frac{A_6}{12} + \frac{a_{111}}{2}, & b_{111} &= \frac{A_8}{6} + \frac{11 A_6}{6} - 2a_{102}, 
\end{aligned}
\end{equation*}
we get 
\begin{equation*}
\begin{aligned}
&  m_{2,9} =\frac{\pi}{60}\left(A_5A_6 + A_5 A_8 - A_6 A_8 + A_7 A_8\right), & m_{2,11} &= \frac{\pi}{90}A_5 A_6 + \frac{8}{3} m_{2,9}, \\
& m_{2,13} = \frac{\pi }{315}A_6 (A_6 - A_7) - \frac{61}{14} m_{2,9} + \frac{26}{7} m_{2,11}, & m_{2,15} &= \frac{365}{56} m_{2,9} - \frac{495}{56} m_{2,11} + 5 m_{2,13}, \\ 
& m_{2,17} = \frac{905}{28} m_{2,9} -\frac{3265}{84} m_{2,11} +\frac{395}{24} m_{2,13}.&& 
\end{aligned}
\end{equation*}
Since $A_5$, $A_6$, $A_7,$ and $A_8$ are free parameters, it is simple to conclude, as above, we get at least $6$ simple zeros of $M_2^{[13]}(r)$, and hence there exists a cubic perturbation of system $CR_1$ exhibiting at least $6$ limit cycles. 
It is here that we disagree with the expansion of $\mathcal{M}_2(r)$ made in the mentioned paper \cite{LiZha2014}, because there the authors got more independence between the terms than we got here. 
	
Anyway, from the expressions in \eqref{4566}, it is clear that if we consider $\alpha$ as a perturbative parameter of $CR_1$ (with $\alpha = 0$), then more order of $\varepsilon$ is needed to thoroughly see its influence, because it appears with exponents up to $12$. We remark here that the center studied in \cite{LiZha2014} corresponds to the one studied here using a convenient rescaling in the $x$ and $y$ variables.

\subsection{The systems $S_1$, $S_2$, $S_3,$ and $S_4$} 
Perturbations of the systems $S_1$, $S_2$, and $S_3$ are common in the literature. 
For $S_4$, although they exist, see for instance \cite{ShiJiaXia2001}, they are not so common. 
Anyway, by following the algorithm of Section~\ref{section:averaging} and similar calculations as the ones done above, it is simple to conclude the raising of a least $2$ limit cycles bifurcating from the origin by using second-order averaging theory for each of these systems. 
Here we sketch these calculations only for $S_4$: 
We assume conditions \eqref{cS4} and calculate $M_2^{[7]}$. 
Then we solve the system $m_{2,2k+1} =A_{k+1}$, $k=0,1,2$, in the coefficients $a_{210}$, $b_{202}$, and $b_{201}$, obtaining 
\[
m_{2,1} = A_1, \quad m_{2,3}  = A_2, \quad m_{2,5}  = A_3, \quad m_{2,7}  = \frac{56}{3} A_3 - \frac{1120}{27} A_2. 
\] 
As above we get at least $2$ limit cycles. 
 
Actually, $2$ is the maximum number of limit cycles of small amplitude that can raise from quadratic perturbations of $S_4$, by a classical result due to Chicone and Jacobs \cite{ChiJac1991}. 

\section{Conclusions and further remarks}
Averaging theory is a powerful method for identifying limit cycles in systems that are nearly integrable, mainly by studying the displacement map. However, in practice, this approach often faces major challenges because the integrals involved are usually complex and hard to calculate, and sometimes impossible to solve. The first-order averaging function can be computed when the level curves of the unperturbed center are properly parametrized; otherwise, we only get integral expressions. The second-order averaging functions are even harder to obtain and can only be expressed as integrals in very rare cases, and only when the first-order function is available in explicit form. Hence, Taylor expansions are usually very useful and, in most cases, the only valid approach.

In this paper, we present a method for handling both first- and second-order averaging functions. The main idea is to find conditions under which the first-order averaging function becomes identically zero, without needing to compute it explicitly. This allows us to move forward and study the second-order averaging function. To do this, we propose an alternative approach that combines the Taylor expansion of the displacement map with the integral form of the averaging function. This strategy helps us establish lower bounds for the cyclicity problem in the context of the Arnold--Hilbert problem. Moreover, by forcing the first-order averaging function to vanish through constraints on the coefficients of the series expansion, we make it possible to carry out a systematic analysis of the second-order averaging function. It is important to note that the number of necessary conditions is not known in advance. This further complicates the computations.

As we have shown, this method has strong potential to simplify the perturbative bifurcation analysis of limit cycles near integrable systems. Still, even though the core idea remains the same, each integrable system brings its own challenges and specific features. Additionally, the method is inherently local and mainly applies to the unfolding of limit cycles that come from first-order averaging.

Despite these limitations, we believe the examples provided, especially the cubic case, demonstrate the strength of the approach. They also help us identify inaccuracies in earlier studies of the $CR_1$ system, as discussed above.

We are confident that the methodology introduced here can offer valuable insights into the study of smooth differential systems. Future work may explore its extension to piecewise-smooth systems and higher-dimensional settings.

\section{acknowledgements} 
We thank the referees for the helpfull comments and suggestions. 
This work was realized under the auspices of the Brazilian S\~ao Paulo Research Foundation FAPESP (grants 2019/07316-0, 2020/14498-4, 2021/14987-8, 2022/14484-9 and 2023/00376-2); the Brazilian agency Conselho Nacional de Desenvolvimento Cient\'ifico e Tecnol\'ogico CNPq (grants 403959/2023-3 and 308112/2023-7); the Brazilian agency Coordena\c{c}\~ao de Aperfei\c{c}oa\-men\-to de Pessoal de N\'ivel Superior CAPES (Finance Code 001); and the Catalan AGAUR agency (grant 2021 SGR 00113); the Spanish AEI agency (grants PID2022-136613NB-I00 and CEX2020-001084-M). 

\section{Conflict of Interest} 
The authors have no conflicts of interest to declare. 


\begin{thebibliography}{10}
	
	\bibitem{Bui2017}
	A.~Buic\u{a}.
	\newblock On the equivalence of the {M}elnikov functions method and the
	averaging method.
	\newblock {\em Qual. Theory Dyn. Syst.}, 16(3):547--560, 2017.
	
	\bibitem{ChaSab1999}
	J.~Chavarriga and M.~Sabatini.
	\newblock A survey of isochronous centers.
	\newblock {\em Qual. Theory Dyn. Syst.}, 1(1):1--70, 1999.
	
	\bibitem{ShiJiaXia2001}
	S.~Chen, J.~Feng, and X.~Zeng.
	\newblock The cyclicity of the period annulus around the quadratic isochronous
	center.
	\newblock {\em Wuhan Univ. J. Nat. Sci.}, 6(4):754--758, 2001.
	
	\bibitem{ChiJac1991}
	C.~Chicone and M.~Jacobs.
	\newblock Bifurcation of limit cycles from quadratic isochrones.
	\newblock {\em J. Differential Equations}, 91(2):268--326, 1991.
	
	\bibitem{Chr2005}
	C.~Christopher.
	\newblock Estimating limit cycle bifurcations from centers.
	\newblock In {\em Differential equations with symbolic computation}, Trends
	Math., pages 23--35. Birkh\"auser, Basel, 2005.
	
	\bibitem{CimGas2020}
	A.~Cima, A.~Gasull, and F.~Ma{\~n}osas.
	\newblock A note on the {Lyapunov} and period constants.
	\newblock {\em Qual. Theory Dyn. Syst.}, 19(1):13, 2020.
	\newblock Id/No 44.
	
	\bibitem{Eca1992}
	J.~\'Ecalle.
	\newblock {\em Introduction aux fonctions analysables et preuve constructive de
		la conjecture de {D}ulac}.
	\newblock Actualit\'es Math\'ematiques. [Current Mathematical Topics]. Hermann,
	Paris, 1992.
	
	\bibitem{FatouP1928}
	P.~Fatou.
	\newblock Sur le mouvement d'un syst{\`e}me soumis {\`a} des forces {\`a}
	courte p{\'e}riode.
	\newblock {\em Bull. Soc. Math. Fr.}, 56:98--139, 1928.
	
	\bibitem{GouTorreGine2021}
	J.~Gin\'{e}, L.~F.~S. Gouveia, and J.~Torregrosa.
	\newblock Lower bounds for the local cyclicity for families of centers.
	\newblock {\em J. Differential Equations}, 275:309--331, 2021.
	
	\bibitem{Ili1998}
	I.~D. Iliev.
	\newblock Perturbations of quadratic centers.
	\newblock {\em Bull. Sci. Math.}, 122(2):107--161, 1998.
	
	\bibitem{Ily1991}
	Y.~S. Ilyashenko.
	\newblock {\em Finiteness theorems for limit cycles}, volume~94 of {\em
		Translations of Mathematical Monographs}.
	\newblock American Mathematical Society, Providence, RI, 1991.
	
	\bibitem{Ily2003}
	Y.~S. Ilyashenko.
	\newblock Centennial history of {H}ilbert's 16th problem.
	\newblock In {\em Fundamental mathematics today}, pages 135--213. Nezavis.
	Mosk. Univ., Moscow, 2003.
	
	\bibitem{LiZha2014}
	S.~Li and Y.~Zhao.
	\newblock Limit cycles of perturbed cubic isochronous center via the second
	order averaging method.
	\newblock {\em Int. J. Bifurcation Chaos Appl. Sci. Eng.}, 24(3):8, 2014.
	\newblock Id/No 1450035.
	
	\bibitem{LliNovTei2014}
	J.~Llibre, D.~D. Novaes, and M.~A. Teixeira.
	\newblock Higher order averaging theory for finding periodic solutions via
	{B}rouwer degree.
	\newblock {\em Nonlinearity}, 27(3):563--583, 2014.
	
	\bibitem{Lou64}
	W.~S. Loud.
	\newblock Behavior of the period of solutions of certain plane autonomous
	systems near centers.
	\newblock {\em Contributions to Differential Equations}, 3:21--36, 1964.
	
	\bibitem{MarRouTon1995}
	P.~Marde\v{s}i\'c, C.~Rousseau, and B.~Toni.
	\newblock Linearization of isochronous centers.
	\newblock {\em J. Differential Equations}, 121(1):67--108, 1995.
	
	\bibitem{Per2001}
	L.~Perko.
	\newblock {\em Differential equations and dynamical systems}, volume~7 of {\em
		Texts in Applied Mathematics}.
	\newblock Springer-Verlag, New York, third edition, 2001.
	
	\bibitem{RomSha2009}
	V.~G. Romanovski and D.~S. Shafer.
	\newblock {\em The center and cyclicity problems. {A} computational algebra
		approach}.
	\newblock Boston, MA: Birkh{\"a}user, 2009.
	
	\bibitem{Rou1998}
	R.~Roussarie.
	\newblock {\em Bifurcation of planar vector fields and {H}ilbert's sixteenth
		problem}, volume 164 of {\em Progress in Mathematics}.
	\newblock Birkh\"{a}user Verlag, Basel, 1998.
	
	\bibitem{Yeu2025}
	M.~Yeung.
	\newblock Dulac's theorem revisited.
	\newblock {\em Qual. Theory Dyn. Syst.}, 24(2):Paper No. 57, 27, 2025.
	
	\bibitem{Zol1994}
	H.~Zoladek.
	\newblock Quadratic systems with center and their perturbations.
	\newblock {\em J. Differential Equations}, 109(2):223--273, 1994.
	
	\bibitem{Zol1995}
	H.~Zoladek.
	\newblock Eleven small limit cycles in a cubic vector field.
	\newblock {\em Nonlinearity}, 8(5):843--860, 1995.
	
\end{thebibliography}

\end{document}